\documentclass[11pt]{amsart}
\usepackage{geometry}                
\usepackage{graphicx,color}
\usepackage{amssymb,amsthm}
\usepackage{epstopdf}

\usepackage{overpic}
\usepackage{rotating}
\usepackage{comment}

\usepackage[colorinlistoftodos,prependcaption,textsize=tiny]{todonotes}

\title[\paa{The Fully} Nonconforming Virtual Element Method for Biharmonic Problems]
{\paa{The Fully} Nonconforming Virtual Element Method \\ for Biharmonic Problems}

\author{P.~F.~Antonietti}
\address{P.~F.~Antonietti: MOX, Dipartimento di Matematica, Politecnico di Milano, Italy}
\email{paola.antonietti@polimi.it}

\author{G.~Manzini}
\address{G.~Manzini: T-5 Group, Theoretical Division Los Alamos National Laboratory, New Mexico, USA; Istituto di Matematica Applicata e Tecnologie Informatiche - CNR, Pavia, Italy}
\email{gmanzini@lanl.gov}

\author{M.~Verani}
\address{M.~Verani: MOX, Dipartimento di Matematica, Politecnico di Milano, Italy}
\email{marco.verani@polimi.it}

\thanks{The authors want to thank Franco Brezzi and Donatella Marini for useful discussions on the topic.
\pa{P. F. Antonietti has been partially supported by the SIR Research Grant no. RBSI14VTOS \emph{PolyPDEs: Non-conforming polyhedral finite element methods for approximation of partial differential equations} funded by MIUR - Italian Ministry of Education, Universities and Research.
G. Manzini has been partially supported by the Laboratory Directed Research and
Development program (LDRD), U.S. Department of Energy Office of Science,
Office of Fusion Energy Sciences, under the auspices of the National
Nuclear Security Administration of the U.S. Department of Energy by
Los Alamos National Laboratory, operated by Los Alamos National Security
LLC under contract DE-AC52-06NA25396.
This work is assigned the LA-UR number LA-UR-16-26955.
M. Verani has been partially supported by the Italian research grant  {\sl Prin 2012}  2012HBLYE4  ``Metodologie innovative nella modellistica differenziale numerica'' and by INdAM-GNCS.
}}
\date{\today}                                           

\newtheorem{lemma}{Lemma}
\newtheorem{theorem}{Theorem}
\newtheorem{remark}{Remark}

\newcommand{\Tau}{\mathcal{T}}
\newcommand{\E}{\mathcal{E}}
\newcommand{\V}{\mathcal{V}}
\newcommand{\N}{\mathcal{N}}

\newcommand{\vrtx}{\mathsf{v}}

\newcommand{\pa}[1]{\textcolor{black}{#1}}
\newcommand{\paa}[1]{\textcolor{black}{#1}}

\begin{document}
\maketitle
\begin{abstract}
In this paper we address the numerical approximation of linear fourth-order elliptic problems on polygonal meshes. In particular, we present a novel nonconforming virtual element discretization of arbitrary order of accuracy for biharmonic problems.  The approximation space is made of {possibly discontinuous functions}, thus giving rise to the {\em fully nonconforming} \paa{virtual element} method. We derive optimal error estimates in \pa{a suitable} (broken) energy norm and present numerical results to assess the validity of the theoretical estimates. 
\end{abstract}

\section{Introduction}
In recent years the study of numerical methods for the approximation of partial differential equations on polygonal and polyhedral meshes has flourished at exponential rate (see, e.g., the special issues of References \pa{\cite{BeiraoErn_2016,Bellomo-Brezzi-Manzini:2014}} for a recent overview of the different methodologies). Among the different proposed methodologies, the Virtual Element Method (VEM),  due to its flexibility in dealing with a wide variety of differential problems, has polarized an increasing research activity. VEM has been introduced in the seminal paper \cite{VEM-basic:2013} and can be seen as an evolution of the Mimetic Finite Difference method, see, e.g., References~\cite{MFD:book, Lipnikov-Manzini-Shashkov:2014} for a detailed description. Since then, the VEM has been proposed to address an increasing number of different problems: general elliptic problems \cite{VEM-elliptic:2016, Berrone-SUPG:2016}, linear and nonlinear elasticity \cite{Vem-elasticity:2013,Paulino-et-al:2014,VEM:inelastic:2015}, plate bending \cite{Brezzi-Marini-plate:2013,Chinosi-Marini:2016}, Cahn-Hilliard \cite{Antonietti-Beirao-Scacchi-Verani:2016}, Stokes \cite{Antonietti-Beirao-Mora-Verani:2014, Beirao-Lovadina-Vacca:2015}, Helmholtz \cite{Perugia-Pietra-Russo:2016}, parabolic \cite{Vacca-Beirao:2015}, Steklov eigenvalue \cite{VEM-Steklov:2015}, elliptic eigenvalue \cite{Gardini-Vacca:2016} and discrete fracture networks \cite{VEM-DFN:2014}. Moreover, several different  variants of the VEM have been developed and analysed: mixed \cite{Brezzi-Falk-Marini:2014,mixed-vem:2016}, discontinuous \cite{VEM-discontinuous}, $H(\text{div})$ and $H(\bold{curl})$-conforming \cite{Hdiv-vem:2016}, hp \cite{Beirao-Chernov-Mascotto-Russo:2016}, serendipity \cite{Serendipity:2016} and nonconforming VEM. This latter formulation has
been first analyzed for elliptic problems~\cite{Ayuso-Lipnikov-Manzini:2016, Cangiani-Manzini-Sutton:2015} and subsequently extended to \pa{the} Stokes problem~\cite{Cangiani-Gyrya-Manzini:2016}. Also, very recently, an approximation method for plate bending problems has been analyzed~\cite{vem-cinesi:2016}, which is based on a globally $C^0$-nonconforming virtual element space. 

In this work we present {\em the fully nonconforming virtual element method} for the approximation of biharmonic 
problems. Our method works on unstructured polygonal meshes, provides arbitrary approximation order and does not require any global $C^0$ regularity for the numerical
solution. 
The numerical approximation of biharmonic problems with nonconforming finite elements on triangular meshes has a very long tradition 
and it is beyond the scope of this introduction to provide a detailed account of it (see, e.g., the classical book \cite{ciarlet:book} for a short \pa{overview}). However, it is worth mentioning that as a by product of the results of this paper we obtain, on triangular meshes, a family of novel nonconforming finite elements of \emph{\pa{arbitrary}} order that are {\em not} continuous. Indeed, for the lowest order our nonconforming virtual element method on simplicial meshes reduces to the classical Morley element \cite{Morley:1968}, while for \pa{higher-order polynomial approximation degrees} it gives rise to a \pa{new} family of nonconforming finite elements.  
  
 \medskip
The outline of the paper is as follows. In Section~\ref{S:1} we recall the continuous problem. In Section~\ref{S:2} we introduce our novel, arbitrary order, nonconforming virtual element discretization for the biharmonic problem. In Section~\ref{S:error} we derive the optimal error estimate in a broken energy norm. In Section~\ref{S:3} we numerically assess the validity of the theoretical estimate and, finally, in Section~\ref{S:4} we draw the conclusions.

\subsection{Notation}
\pa{Throughout the paper we shall use the standard notation of the Sobolev spaces $H^m(\mathcal{D})$ for 
a nonnegative integer $m$ and an open bounded domain $\mathcal{D}$. The $m$-th seminorm 
of the function $v$ will be defined by
$$ \vert v \vert_{m,\mathcal{D}} = \sum_{\vert \alpha \vert = m} 
\bigg\| \frac{\partial^{\vert \alpha\vert} v}{\partial_{x_1}^{\alpha_1} \partial_{x_2}^{\alpha_2}}\bigg\|^2_{0,\mathcal{D}},$$
where 
$\| \cdot \|^2_{0,\mathcal{D}}$ stands for the $L^2(\mathcal{D})$ norm and we set
$\vert {\alpha}\vert = \alpha_1 +\alpha_2$ 
for the nonnegative multi-index ${\alpha} =(\alpha_1,\alpha_2)$.
For any integer $m\geq 0$, $\mathbb{P}^m ({\mathcal D})$ is the classical space of polynomials of total degree up to $m$ defined on ${\mathcal{D}}$.
Moreover, $n=(n_1,n_2)$ is the outward unit normal vector to $\partial \mathcal{D}$, the boundary of $\mathcal{D}$, and $t=(t_1,t_2)$ the unit tangent vector in the counterclockwise orientation of the boundary.
To ease the notation, we may use $u_{,i}$ to indicate the first order derivative along the $i$-th direction, and, accordingly, $u_{,n}$ and $u_{,t}$ for the normal and tangential derivatives. Whenever convenient, we shall also use the notation $\partial_n u$ and $\partial_t u$ instead of $u_{,n}$ and $u_{,t}$.
Moreover, we may denote high-order derivatives by repeating the index subscripts, e.g.,
$u_{,ij}={\partial^2 u}\slash{\partial x_i\partial x_j}$, and, likewise, $u_{,nn}$, $u_{,tt}$, $u_{,nnt}$, etc, for multiple derivatives in the normal and tangential directions.
We also use the summation convention of repeated indexes (Einstein's convention), so that
$$u_{,ij}v_{,ij}=\sum_{ij=1}^{2}\frac{\partial^2 u}{\partial x_i\partial x_j}\frac{\partial^2 v}{\partial x_i\partial x_j}.$$
Finally, the notation $A\lesssim B$ will signify that $A\leq c B$ for some positive constant $c$ independent of the discretization parameters.
 \\}

\medskip

\section{The continuous problem}\label{S:1}

Let $\Omega\subset \mathbb{R}^2$ be a convex polygonal domain occupied by the plate with boundary $\Gamma$ and let $f\in L^2(\Omega)$ be a transversal load acting on the plate. According to the Kirchoff-Love model for thin plates~\cite{Landau-Lifshitz:book} and assuming that the plate is clamped all over the boundary, the transversal displacement $u$ is solution to the following problem
\begin{subequations}\label{pb:continuous}
\begin{eqnarray}
D\Delta^2 u &=& f \qquad \text{in~}\Omega\label{pb:1}\\
u&=&0 \qquad \text{on~}\Gamma\label{pb:2}\\
\partial_n u&=&0 \qquad \text{on~}\Gamma\label{pb:3}
\end{eqnarray}
\end{subequations}
where $D=\frac{Et^3}{12(1-\nu^2)}$ is the bending rigidity, $t$ being the thickness, $E$ the Young modulus, and $\nu$ the Poisson's ratio. 

Consider the functional space $V=\big\{v\in H^2(\Omega):v=\partial_n v=0\text{~on~}\Gamma\big\}$ \paa{and 
denote by $\langle \cdot, \cdot \rangle$ the duality pairing between $V$ and its dual $V^*$.}

The variational formulation of \eqref{pb:continuous} reads as: \emph{Find $u\in V$ such that}
\begin{equation}\label{pb:w}
a(u,v) = \paa{\langle f, v \rangle} \qquad\forall v\in V,
\end{equation}
where 
\begin{align}
a(u,v)=D \int_\Omega \big(\nu \Delta u\Delta v + (1-\nu)\,u_{,ij} v_{,ij}\big) dx 
\quad\textrm{and}\quad 
\paa{\langle f, v \rangle}=\int_\Omega f v dx.
\label{pb:bil-rhs}
\end{align}


Setting $\| \cdot \|_V=\vert \cdot \vert_{2,\Omega}$ and employing the boundary conditions and Poincar\'e inequality, 
we can prove that  $\| \cdot \|_V$ is a norm on $V$. Moreover, it holds that 
\begin{subequations}
\begin{align}
a(v,v)            &\gtrsim  \|v\|_{V}^2\phantom{\|u\|_{V}\|v\|_{V}}\forall v\in V \label{coerc}\\[0.25em]
\vert a(u,v)\vert &\lesssim \|u\|_{V}\|v\|_{V}\phantom{\|v\|_{V}^2}\forall u,v\in V. \label{cont}
\end{align}
\end{subequations}
Hence, there exists a unique solution $u\in V$ to \eqref{pb:w} (see, e.g., \cite{Brenner-scott:book}).

\subsection{\pa{Preliminaries}}

In this section, we collect some useful definitions that will be employed in the rest of the paper. Let $\sigma_{ij}(u)=\lambda (u_{,11}+u_{,22})\delta_{ij}+ \mu u_{,ij}$ with Lam\'e parameters $\lambda=D\nu$, $\mu=D(1-\nu)$. We set 
\begin{equation}\label{T:1}
M_{nn}(u)= \sigma_{ij}\paa{(u)} n_i n_j,\quad
M_{nt}(u)=\sigma_{ij} \paa{(u)}  n_i t_j,\quad
T(u)=\sigma_{ij,j} \paa{(u)} n_i+M_{nt,t} \paa{(u)} ,
\end{equation}
(we recall the summation notation of repeated indexes)
and observe that
\begin{equation}\label{T:6}
\begin{aligned}
M_{nn}(u)&=\Delta u - (1-\nu) u_{,tt} = \nu \Delta u + (1-\nu) u_{,nn},\\
M_{nt}(u)&=u_{,nt},\\
T(u) &= \partial_n(\Delta u) + (1-\nu) u_{,ntt}.
\end{aligned}
\end{equation}
Moreover, 
let $K\subset \mathbb{R}^2$ be a polygonal domain and set 
$$a^K(u,v)=D \int_K \big(\nu \Delta u\Delta v + (1-\nu)\,u_{,ij} v_{,ij}\big)dx.$$ 
Integrating by parts and employing \eqref{T:1} and \eqref{T:6} yield the following useful identities

\begin{align}\label{eq:1} 
a^K(u,v) 
= D\bigg\{ &
\int_K \Delta^2 u v dx + \int_{\partial K} \big(\Delta u - (1-\nu) u_{,tt}\big) v_{,n} ds \nonumber \\[0.5em]
&- \int_{\partial K} \big(\partial_n(\Delta u) v - (1-\nu) u_{,nt} v_{,t} \big) ds 
\bigg\}\nonumber\\[1em]
= D\bigg\{ &\int_K \Delta^2 u v dx  + 
\int_{\partial K} M_{nn}(u) \partial_n v \,ds - \int_{\partial K}T(u)v ds - \sum_{e\in \partial K} (M_{nt}(u),v n_{\partial e})_{\partial e}\bigg\}
\end{align} 
where $\partial e$ is the boundary of edge $e\subseteq \partial K$ and $n_{\partial e}$ is the outwards normal ``vector'' to $\partial e$. 
\pa{For every edge $e$ with end points $\vrtx_1$ and $\vrtx_2$ the boundary ${\partial e}$ is the set $\{\vrtx_1,\vrtx_2\}$, and, depending on the chosen edge orientation, $n_{\partial e}$ at the end points is equal to $+1$ or $-1$. }

\section{Nonconforming virtual element discretization} \label{S:2}

\pa{
The nonconforming virtual element approximation of the variational problem~\eqref{pb:continuous} reads as: \emph{Find $u_h\in V_{h,\ell}$ such that} 
\begin{equation}\label{pb:ncVEM}
a_h(u_h,v_h)=\langle f_h,v_h\rangle
\qquad \forall v_h\in V_{h,\ell},
\end{equation}
where $V_{h,\ell}$ is the nonconforming virtual element space 
of order $\ell$ that approximates the functional space $V$, and $a_h(\cdot,\cdot)$ 
and $\langle f_h, \cdot\rangle$ are the nonconforming virtual element 
bilinear form and load term that approximate $a(\cdot,\cdot)$ and 
$\langle f, \cdot \rangle$ in~\eqref{pb:w}, respectively.   
The rest of this section is devoted to the construction of these quantities.}

\subsection{Technicalities}

Let $\{\Tau_h\}_{h}$ be a sequence of decompositions (meshes) of $\Omega$ into non-overlapping polygons $K$.
Each mesh $\Tau_h$ is labeled by the mesh size parameter $h$, which will be defined below, and satisfies a few regularity assumptions that are necessary to prove the convergence of the method and derive an estimate of the approximation error. These regularity assumptions are introduced and discussed in Section~\ref{S:error}. 
Let $\E_h$ be the set of edges in $\Tau_h$ such that $\E_h=\E_h^i \cup \E_h^{\Gamma}$, where $\E_h^i$ and $\E_h^{\Gamma}$ are the set of interior and boundary edges, respectively. Similarly, we denote by $\V_h=\V_h^i \cup \V_h^{\Gamma}$ the set of vertices in $\Tau_h$, where $\V_h^i$ and $\V_h^{\Gamma}$ are the sets of interior and boundary vertices, respectively. Accordingly, $\V_h^K$ is the set of vertices of $K$. Moreover, $|K|$ and $|e|$ denotes the area of cell $K$ and the length of edge $e$, $\partial K$ is the boundary of $K$, $h_K$ is the diameter of $K$ and the mesh size parameter is defined as $h=\max_{K\in\Tau_h}h_K$.  
We introduce the broken Sobolev space for any integer number $s>0$
$$ H^s(\Tau_h)=\Pi_{K\in \Tau_h} H^s(K) =\big\{ v\in L^2(\Omega): v_{\vert K} \in H^s(K) \textrm{ for any } K \in \Tau_h\big\}$$
\pa{and endow it with the broken $H^s$-seminorm
 $ |v|^2_{s,h}=\sum_{K\in \Tau_h} | v|^2_{s,K}$.
}
%
\pa{We denote the traces of $v$ 
on $e\subset \partial K^+\cap \partial K^-$ 
from the interior of $K^\pm$ by $v^\pm$, respectively. Then, 
we define the jump of $v$ on the interior edge $e \in \E_h^i$ ƒby
$[v]=v^+ - v^-$ and on the boundary edge $e\in\mathcal{E}_h^{\Gamma}$ by $[v]=v_{|e}$.}

For future use, we also introduce the nonconforming space $H^{2,nc}(\Tau_h)\subset H^2(\Tau_h)$ defined as follows
\begin{align}
H^{2,nc}(\Tau_h)=\bigg\{ &v \in H^2(\Tau_h): v \text{~continuous~at~internal~vertexes}, 
\,\, v_h(\vrtx_i)=0  ~\forall \vrtx_i \in \mathcal{V}_h^{\Gamma}\nonumber\\
&\int_e [\partial_n v] ds =0~ \forall e\in \mathcal{E}_h \bigg\}. \nonumber
\end{align}
\pa{We next prove the following result.
\begin{lemma}\label{eq:norm2broken}
$\vert \cdot \vert_{2,h}$ is a norm on both $V$ and $H^{2,nc}(\Tau_h)$.
\end{lemma}
\begin{proof}
Employing \cite[Corollary 4.2]{Brenner-H2-broken:2004} and 
\cite[(5.2)]{Brenner-H2-broken:2004} (with $\Phi(v)$ chosen as in 
\cite[Example 2.6]{Brenner-H2-broken:2004}) we can prove that 
\begin{equation*}\label{eq:poincare}
\vert v\vert_{1,h}\lesssim \vert v \vert_{2,h}\qquad \forall v\in H^{2,nc}(\Tau_h)
\end{equation*}
which implies that $\vert v\vert_{2,h}$ is a norm on $H^{2,nc}(\Tau_h)$. 
\end{proof}}
\subsection{Local and global nonconforming virtual element space}
In this section, we introduce the local and global nonconforming virtual element spaces. \\

\noindent
For $\ell\geq 2$, the local virtual element space is defined as follows:
\begin{equation*}
V_{h,\ell}^K=\big\{v_h \in H^2(K): \Delta^2 v_h \in \mathbb{P}^{\ell-4}(K), M_{nn}(v_h)_{|e}\in \mathbb{P}^{\ell-2}(e), T(v_h)\in \mathbb{P}^{\ell-3}(e)~\forall e \in \partial K
\big\}
\end{equation*}
with the usual convention that $\mathbb{P}^{-1}(K)=\mathbb{P}^{-2}(K)=\{ 0 \}$. The solution of the biharmonic problem in the definition of $V_{h,\ell}^K$ is uniquely determined up to a linear function that can be filtered out by fixing the value at three non-aligned vertexes of $K$. 

\begin{remark}
By construction, it holds that $\mathbb{P}^\ell(K) \subset V_{h,\ell}^K$.
\end{remark}

\begin{figure}
\centering
\begin{tabular}{cccc}
\includegraphics[scale=0.2]{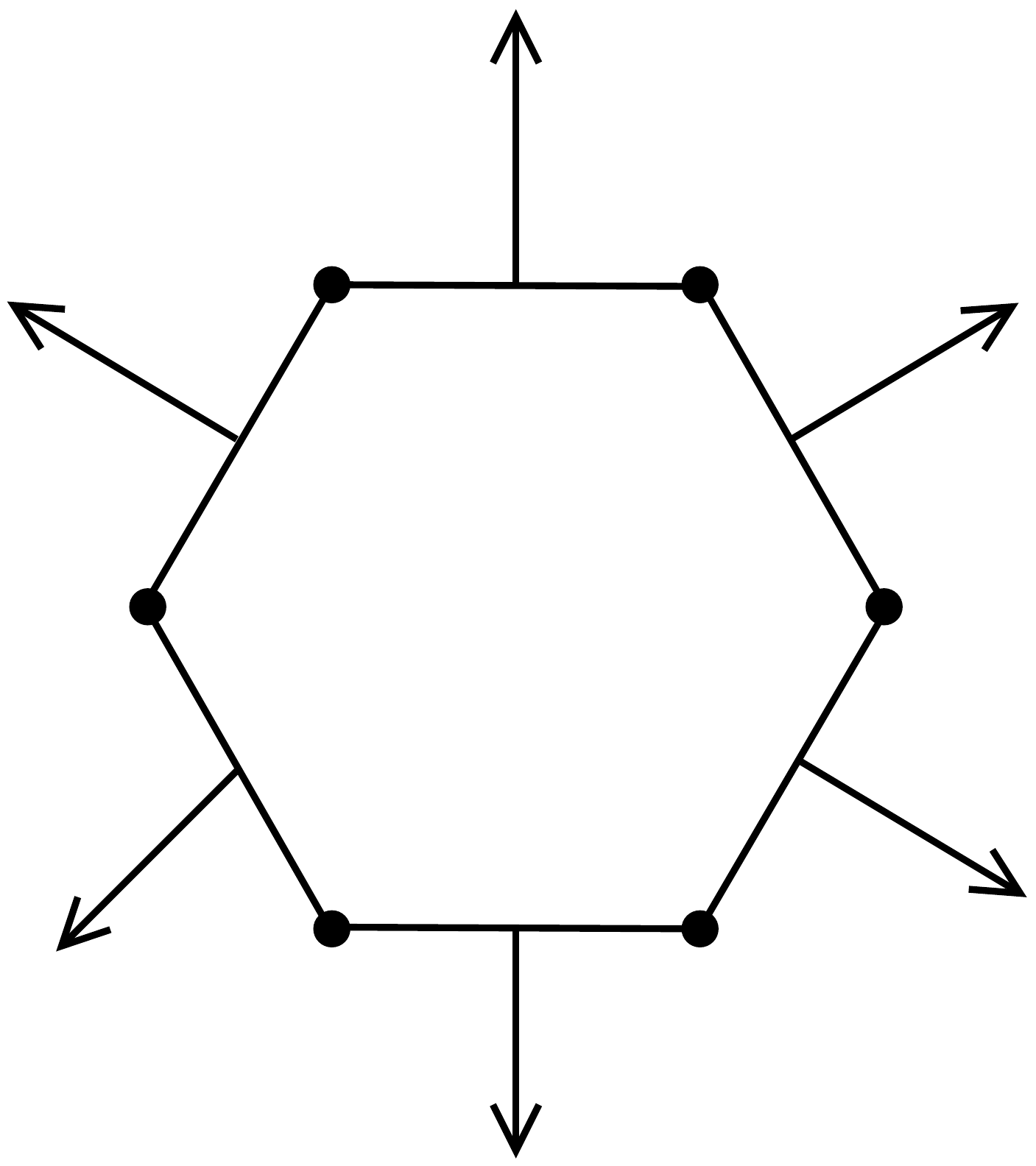} &
\includegraphics[scale=0.2]{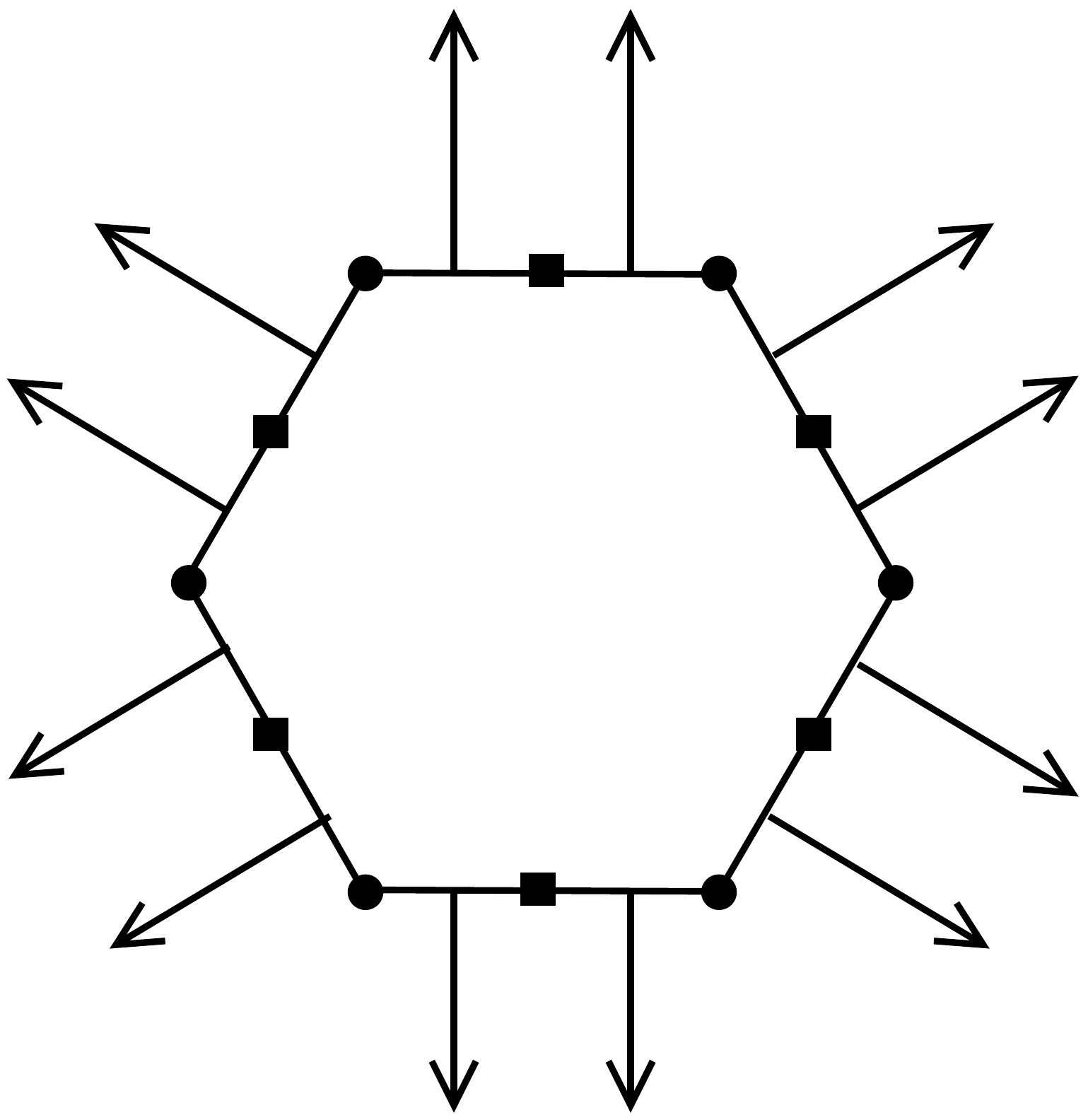} &
\includegraphics[scale=0.2]{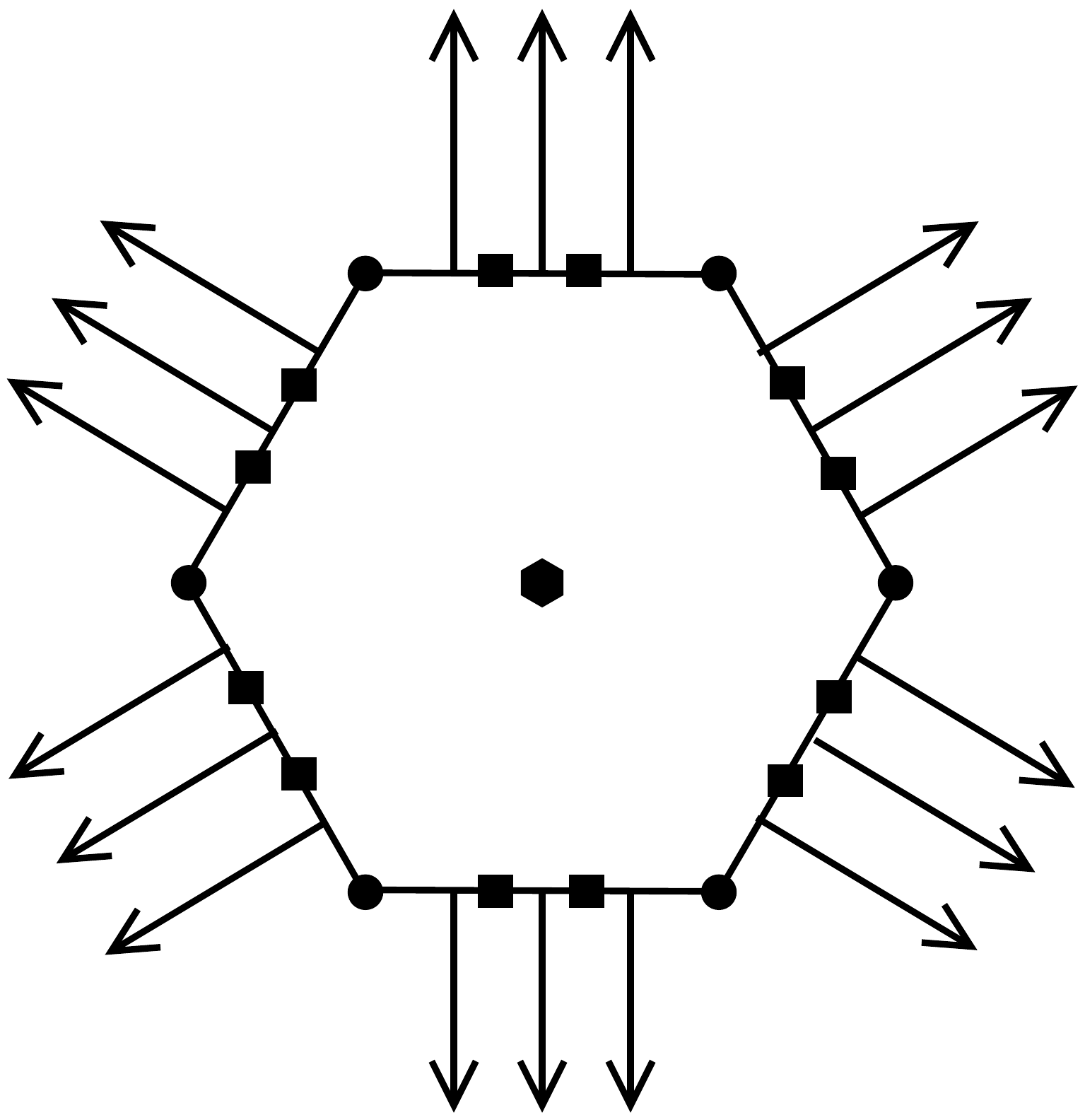} &
\includegraphics[scale=0.2]{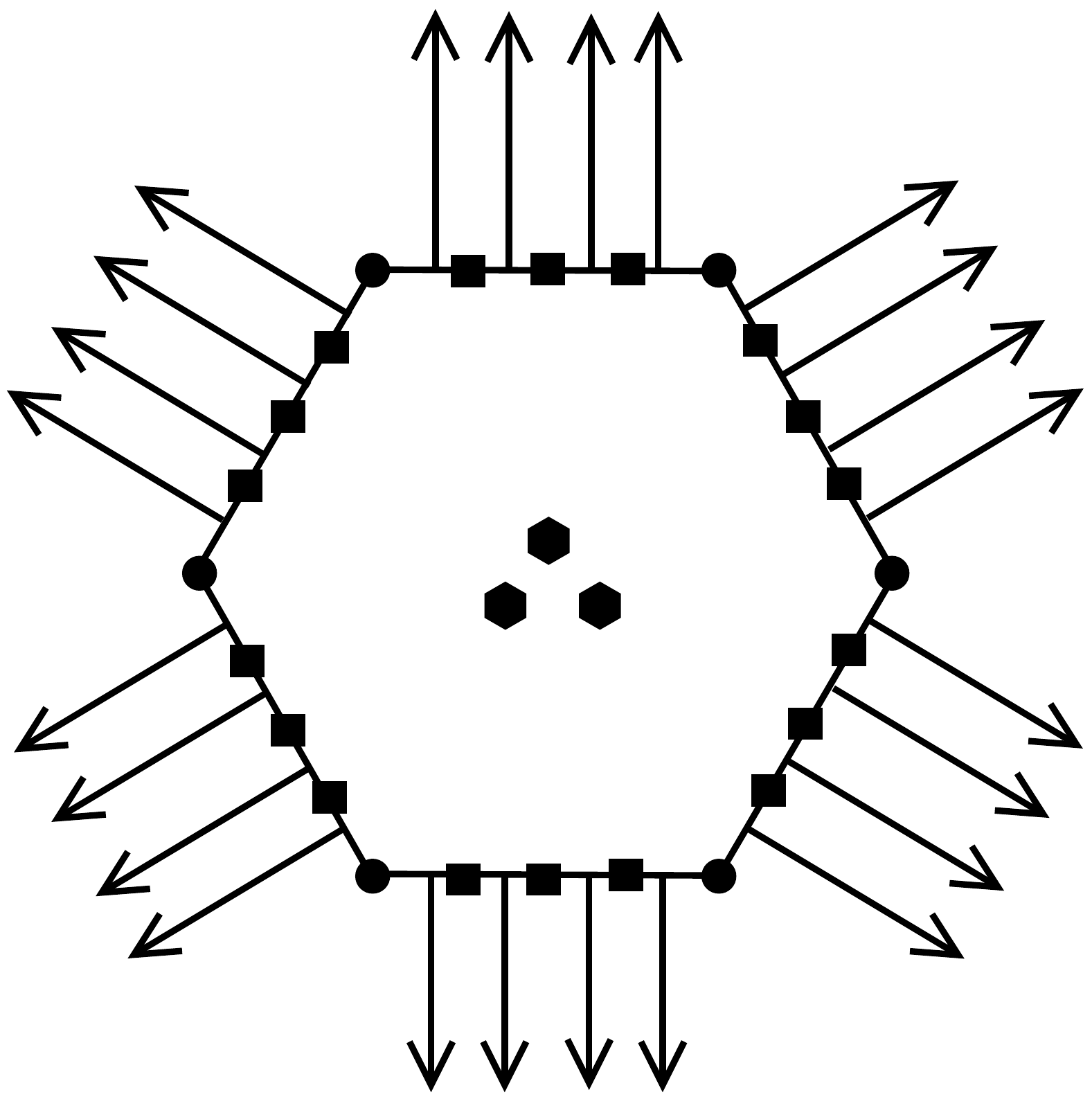} \\[0.2em]
\textbf{$\ell=2$} & \textbf{$\ell=3$} & \textbf{$\ell=4$} & \textbf{$\ell=5$}
\end{tabular}
\caption{Local degrees of freedom of $V^K_{h,\ell}$ for $l=2,3,4,5$: vertex values 
$u_h(\vrtx)$ (dots); edge moments of $\partial_n u_h$ (arrows); edge moments of $u_h$ 
(squares); cell moments of $u_h$ (central dots).}
\label{fig:dofs:2-5}
\end{figure}

\noindent
We choose the \emph{degrees of freedom} of $V_{h,\ell}^K$ as follows:

\smallskip
\begin{itemize}
\item[(D1)] for $\ell\geq 2$: $v_h(\vrtx_i)$ for any vertex $\vrtx_i$ of $K$;

\medskip
\item[(D2)] for $\ell\geq 2$: $\displaystyle\int_e \partial_n v_h\, p\, ds$ for any $p\in \mathbb{P}^{\ell-2}(e)$ and any edge $e$ of $\partial K$;

\medskip
\item[(D3)] for $\ell\geq 3$: $\displaystyle{\frac{1}{\vert e\vert}}\int_e v_h p\, ds$ for any $p\in \mathbb{P}^{\ell-3}(e)$ and any edge $e$ of $\partial K$;

\medskip
\item[(D4)] for $\ell\geq 4$: $\displaystyle{\frac{1}{\vert K\vert}}\int_K v_h p\, ds$ for any $p\in \mathbb{P}^{\ell-4}(K)$.
\end{itemize}

\pa{These degrees of freedom are illustrated in Figure~\ref{fig:dofs:2-5} for 
the virtual element spaces with $\ell=2,3,4,5$. We next show that these degrees of freedom are unisolvent in $V^{K}_{h,\ell}$.}
\begin{lemma}\label{lemma:unisolvence}
The degrees of freedom (D1)-(D4) are unisolvent for  $V_{h,\ell}^K$.
\end{lemma}
\begin{proof}
Employing \eqref{eq:1}, for any $v_h\in V_{h,\ell}^K$ there holds 
\begin{align}
a^K(v_h,v_h) 
= D\bigg\{ &\int_K \Delta^2 v_h \,v_h dx
+\int_{\partial K} M_{nn}(v_h) \partial_n v_h ds  
-\int_{\partial K} T(v_h) \partial_n v_h ds \nonumber\\
& -\sum_{e\in \partial K} (M_{nt}(v_h),v_h n_{\partial e})_{\partial e}
\bigg\}.\label{aux2:unisolv}
\end{align}
We first observe that $\Delta^2 v_h$ is a polynomial of order $\ell-4$ on $K$. Moreover, we note that on each edge the functions $M_{nn}(v_h)$ and $T(v_h)$ are polynomial of degree up to $\ell-2$ and $\ell-3$, respectively. Hence, by setting to zero the degrees of freedom we get   $a^K(v_h,v_h)=0$.
This latter implies $\vert v_h\vert_{H^2(K)}=0$ which gives $v_h=0$ in view of the fact that setting to zero the nodal values of $v_h$ filters the linear polynomials (i.e. the kernel of the seminorm $\vert \cdot\vert_{H^2(K)}$). 
\end{proof}


Building upon the local spaces $V_{h,\ell}^K$, the global nonconforming virtual element space is then defined as follows
\begin{align}
V_{h,\ell}=\bigg\{ & v_h: v_h\vert_K\in V_{h,\ell}^K, v_h \text{~continuous~at~internal~vertexes},v_h(\vrtx_i)=0  ~\forall \vrtx_i \in \mathcal{V}_h^{\Gamma} \nonumber\\
& \int_e [\partial_n v_h] p ds =0~ \forall p\in \mathbb{P}^{\ell-2}(e), 
\int_e [v_h] p ds =0 ~ \forall p\in \mathbb{P}^{\ell-3}(e) ~\forall e\in \mathcal{E}_h \bigg\}. \label{global-space}
\end{align}
We observe that by construction it holds $V_{h,\ell}\subset H^{2,nc}(\Tau_h)$ and $V_{h,\ell}\not \subseteq H^2_0(\Omega)$. Moreover, it is important to remark that
our nonconforming virtual element space does not require that its functions are globally continuous over $\Omega$, thus admitting piecewise discontinuous functions on each partition $\Tau_h$ (see also Remark \ref{rem:1} below). 

\begin{remark}[Lowest order case]\label{rem:1} Let us briefly comment on the lowest-order VE space, for $\ell=2$. 
In this case the local space is given by
\begin{equation*}
V_{h,2}^K=\big\{ v_h \in H^2(K): \Delta^2 v_h =0, M_{nn}(v_h)_{|e}\in \mathbb{P}^0(e), T(v_h)_{|e}=0\quad \forall e \in \partial K\big\},
\end{equation*}
and the local degrees of freedom are: 
\smallskip
\begin{itemize}
\item[\emph{(d1)}] $v_h(\vrtx_i)$ for any vertex $\vrtx_i$ of $K$;

\smallskip
\item[\emph{(d2)}]$\displaystyle\int_e \partial_n v_h ds$ for any edge $e$ of $\partial K$.
\end{itemize} 
On triangular meshes the degrees of freedom (d1)-(d2) of $V_h^K$ are the same of the Morley's nonconforming finite element space ~\cite{Morley:1968} and the unisolvence property from Lemma~\ref{lemma:unisolvence} implies that these two local spaces coincide.  
%
Finally, the global lowest-order nonconforming virtual element space is given by
\begin{align*}
V_{h,2}=\bigg\{&v_h: v_h\vert_K\in V_h^K, \,v_h \text{~continuous~at~internal~vertexes},
\,v_h(\vrtx_i)=0 ~\forall \vrtx_i \in \mathcal{V}_h^{\Gamma},\\ 
&\textrm {and}\int_e [\partial_n v_h] ds =0~ \forall e\in \mathcal{E}_h\bigg\}
\end{align*}
and clearly contains functions that are piecewise discontinuous on $\Tau_h$.
\end{remark}

\subsection{Construction of the bilinear form}

\paa{Starting from the local
bilinear forms $a_h^K(\cdot,\cdot)\,:\,V^K_{h,\ell}\times V^K_{h,\ell}\to \mathbb{R}$ the global bilinear form $a_h(\cdot,\cdot)$ is assembled in the usual way:}
$$ a_h(u_h,v_h) =\sum_{K\in \Tau_h} a_h^K(u_h,v_h).$$
Each local bilinear form is given by 
\begin{equation}\label{local-bilinear}
a_h^K(u_h,v_h)= a^K\Big(\Pi ^{\Delta,K}_{\ell} u_h, \Pi ^{\Delta,K}_{\ell} v_h\Big) + S^K \Big(\big(I- \Pi ^{\Delta,K}_{\ell}\big) u_h, \big(I- \Pi ^{\Delta,K}_{\ell}\big) v_h \Big),
\end{equation}
where $\Pi ^{\Delta,K}_{\ell}$ is the elliptic projection operator discussed below and $S^K(u_h,v_h)$ is a symmetric and positive definite bilinear form such that
$$ a^K(v_h,v_h) \lesssim S^K(v_h,v_h) \lesssim a^K(v_h,v_h) $$
for all $ v_h\in V_{h,\ell}^K$ such that $\Pi ^{\Delta,K}_{\ell} v_h=0$. 
A practical and very simple choice for $S^K(\cdot ,\cdot)$ is the Euclidean scalar product associated to the degrees of freedom scaled by factor $h_k^{-2}$.

\medskip
The operator $\Pi^{\Delta,K}_{\ell} \colon V^K_{h,\ell} \rightarrow \mathbb{P}^\ell(K)$ is the solution of the elliptic projection problem:
\begin{align}
a^K(\Pi ^{\Delta,K}_{\ell} v_h, p) 
&= a^K(v_h, p)\phantom{(\!( v_h, p)\!)_K}\forall p\in \mathbb{P}^\ell(K),\label{def:proj:1}\\[0.5em]
(\!( \Pi ^{\Delta,K}_{\ell} v_h , p )\!)_K 
&= (\!( v_h, p)\!)_K\phantom{a^K(v_h, p)}\forall p\in \mathbb{P}^1(K),\label{def:proj:2}
\end{align}
where 
$$
(\!( v_h, w_h )\!)_K {=} \sum_{\vrtx\in\V^K_h} \!\! v_h(\vrtx) \: w_h(\vrtx).
$$
It is immediate to verify that $\Pi ^{\Delta,K}_{\ell}$ is a projector onto the space of polynomials $\mathbb{P}^\ell(K)$. Indeed, for any $q\in \mathbb{P}^{\ell}(K)$ equation \eqref{def:proj:1} with   $p=\Pi^{\Delta,K}_{\ell}q-q$ yields $(\Pi^{\Delta,K}_{\ell} q)_{,ij} = q_{,ij}$ for $i,j=1,2$. This latter relation combined with \eqref{def:proj:2} proves the assertion. Furthermore, as stated by the following lemma, the polynomial projection $\Pi ^{\Delta,K}_{\ell}v_h$ is computable from the degrees of freedom of $v_h$.

\begin{lemma}
The projector $\Pi ^{\Delta,K}_{\ell} \colon V^K_{h,\ell} \rightarrow \mathbb{P}^\ell(K)$ can be computed using only the degrees of freedom (D1)-(D4). 
\end{lemma}
\begin{proof} 
In view of \eqref{def:proj:1} and assuming, as usual, the computability of $a^K(p,q)$ for polynomial functions $p,q$, it is sufficient to prove the computability of $a^K(p,v_h)$ for any $p\in \mathbb{P}^\ell(K)$ and $v_h\in V_{h,\ell}^K$.
%
%
%
%
Employing \eqref{eq:1} we have 
\begin{align*}
a^K(p,v_h) =  
D\bigg\{ &\int_K \Delta^2 p v_h dx 
+ \int_{\partial K} M_{nn}(p) \partial_n v_h ds 
-\int_{\partial K} T(p) v_h ds\\[0.5em]
& -\sum_{e\in \partial K} \big(M_{nt}(p),v_h n_{\partial e}\big)_{\partial e} 
\bigg\}.
\end{align*} 
Each term of the right-hand side can be computed using only the degrees of freedom (D1)-(D4). 
Indeed, for the first term we note that $\Delta^2 p$ is a polynomial of order $\ell-4$; for the second term we note that $M_{nn}(p)$ is a polynomial of order $\ell-2$ on each edge; 
for the third term we note that $T(p)$ is a polynomial of degree $\ell-3$; finally, we note that the last term 
depends on the value of $v_h$ at the vertexes of $K$.
\end{proof}

The local bilinear form $a_h^K$ has the two crucial properties of polynomial consistency and 
stability that we state in the following lemma. 

\medskip
\begin{lemma}~\\
\label{lemma:properties}
\vspace{-0.5\baselineskip}
\begin{itemize}
\item {\bf $\ell$-consistency}: For any $ p\in\mathbb{P}^\ell(K)$ and any 
$ v_h\in V^K_{h,\ell}$ it holds that: 
\begin{equation}\label{eq:consistency}
a_h^K(p,v_h)=a^K(p,v_h).
\end{equation}

\medskip
\item {\bf stability}: For any $ v_h\in V^K_{h,\ell}$ it holds that:
\begin{equation}\label{eq:stability}
a^K(v_h,v_h) \lesssim a_h^K(v_h,v_h) \lesssim a^K(v_h,v_h),
\end{equation}
where the hidden constants are independent of $h$ and $K$ (but may depend on $\ell$). 
\end{itemize}
\end{lemma}

The proof is straightforward, and is therefore omitted.

\subsection{Construction of the load term}

%
%
Let $\Pi_K^{\ell}$ denote the $L^2$-projection onto $\mathbb{P}^{\ell}(K)$ and $f_h$ 
be the piecewise polynomial approximation of $f$ on $\Tau_h$ given by 
\begin{equation}\label{eq:rhs}
{f_h}_{|K}=\Pi_K^{\ell-2} f
\end{equation}
for $\ell\geq 2$ and every $K \in \Tau_h$. Then, we set 
\begin{equation}\label{vem:rhs}
\langle f_h,v_h\rangle=\sum_{K\in \Tau_h} \int_K f_h v_h\,dx.
\end{equation}
In view of \eqref{eq:rhs} and using the definition of the $L^2$-projection we find that
\begin{equation}\label{aux:1}
\langle f_h,v_h\rangle=\sum_{K\in \Tau_h} \int_K  \Pi_K^{\ell-2} f v_h\,dx=\sum_{K\in \Tau_h} \int_K \Pi_K^{\ell-2} f \Pi_K^{\ell-2} v_h\,dx=
\sum_{K\in \Tau_h} \int_K f \Pi_K^{\ell-2} v_h\,dx. 
\end{equation}
The right-hand side of \eqref{aux:1} is computable by using the degrees of freedom (D1)-(D4) and the enhanced approach~\cite{Enhanced-VEM:2013} 
that considers the augmented local space 
\begin{align*}
W_{h,\ell}^K =\bigg\{ & v_h \in H^2(K): \Delta^2 v_h \in \mathbb{P}^{\ell-2}(K), \,M_{nn}(v_h)_{|e}\in \mathbb{P}^{\ell-2}(e), \,T(v_h)\in \mathbb{P}^{\ell-3}(e)~\forall e \in \partial K,\\
&\int_K \Pi ^{\Delta,K}_{\ell} v_h p dx = \int_K v_h p dx ~~\forall p\in \mathbb{P}^{\ell-2}\setminus \mathbb{P}^{\ell-4}
 \bigg\}.
 \end{align*}
Since (D1)-(D4) are still unisolvent in $W_{h,\ell}^K$, we can compute the projection 
$\Pi^{\ell-2}_{K}$ from the degrees of freedom of $v_h$. 

Finally, from \eqref{aux:1}, employing the Cauchy-Schwarz inequality, standard approximation error estimates and \eqref{eq:norm2broken} we have the estimate
\begin{equation}\label{rhs:error}
\langle f-f_h,v_h\rangle=\sum_{K\in \Tau_h} 
\int_K \big(I-\Pi_K^{\ell-2}\big)f\,\,\big(I-\Pi_K^{0}\big) v_h \,dx 
\lesssim  h^{\ell} |v_h|_{2,h},
\end{equation}
which will be useful in the error analysis of the next section.
  
\section{\pa{Error estimates}}\label{S:error}

We now turn to the derivation of an optimal error estimate for the virtual element discretization \eqref{pb:ncVEM}. \\

\pa{
On the mesh sequence $\{\Tau_h\}_h$ we make the following regularity 
assumptions:
\begin{itemize}
\item [(\bf{H})] there exists a fixed number $\rho_0>0$ independent of $\Tau_h$, such that for every element $K$ it holds:
\smallskip
\begin{enumerate}
\item[(H1)] $K$ is star-shaped with respect to all the points of a ball of radius $\rho_0 h_K$
\smallskip
\item[(H2)] every edge $e\in \E_h$ has length $\vert e \vert \geq \rho_0 h_K$.
\smallskip
\item[(H3)] \paa{There exists a point $x_B$ interior to $K$ such that the sub-triangulation obtained by connecting $x_B$ to the vertices of $K$ is made of shape regular triangles.}
\end{enumerate}
\end{itemize}
}
\paa{The assumptions (H1)-(H2) are standard (see, e.g., ~\cite{VEM-basic:2013}) while (H3) is required to perform the error analysis (see, in particular, \eqref{H3-eq-1} and \eqref{H3-eq-2} in the proof of Theorem \ref{th:conv}).}

\pa{In view of \paa{the assumptions (H1)-(H2)} on $\Tau_h$, we can define, for every
smooth enough function $w$, an ``interpolant'' in $V_{h,\ell}$ with the right interpolation properties.
More precisely, if $\chi_i(w)$, $i = 1,\ldots,\mathcal{N}$, denotes the $i$-th global degree of freedom of a sufficiently regular function $w$, there exists a
unique element $w^I \in V_{h,\ell}$ such that
$$\chi_i(w-w^I) = 0 \qquad i = 1,2,\ldots, \mathcal{N}.$$
Moreover, combining Bramble-Hilbert technique and scaling
arguments (see e.g. \cite{VEM-basic:2013,VEM-Steklov:2015} and \cite{Brenner-scott:book}) as in the  
finite element framework we can prove that
\begin{equation*}
\|w-w^I\|_{s,\Omega}\lesssim C h^{\beta - s} |w |_{\beta,\Omega} \qquad s=0,1,2 \quad 3\leq \beta \leq k+1.
\end{equation*}
}
In accordance with the seminal paper \cite{lascaux-lesaint:1975} (see also \cite{ciarlet:book}) we obtain the following result.

\begin{theorem}
Under the regularity mesh assumptions \paa{(H1)-(H2)}, there exists a unique solution $u_h\in V_{h,\ell}$ to \eqref{pb:ncVEM}. Moreover, for every approximation $u_\pi \in \mathbb{P}^\ell(\Tau_h)$ of the exact solution $u$ of \eqref{pb:w}, it holds that
\begin{equation}\label{vem-cea}
| u-u_h |_{2,h}\lesssim (| u-u^I |_{2,h} + | u-u_\pi |_{2,h} + \sup_{v_h \in V_{h,\ell}} \frac{\langle f-f_h,v_h\rangle}{| v_h |_{2,h}} + \sup_{v_h\in V_{h,\ell}}
\frac{\N(u,v_h)}{| v_h|_{2,h}},
\end{equation}
where $u^I\in V_{h,\ell}$ is the interpolant of $u$ in the virtual element space $V_{h,\ell}$ and
\begin{align}
\N(u,v_h) 
&= \langle f, v_h \rangle - \sum_{K \in \Tau_h} a^K (u,v_h)
\nonumber\\[0.5em]
&
= D\sum_{K\in \Tau_h}\bigg\{ \int_{\partial K} (\Delta u - (1-\nu) u_{,tt}) {v_h}_{,n} ds 
  -\int_{\partial K} \left(\partial_n(\Delta u) v_h - (1-\nu) u_{,nt} {v_h}_{,t} \right)ds \bigg\}
\label{def:N}
\end{align}
is the non-conformity error.
\end{theorem}
\begin{proof}
Existence and uniqueness of the discrete solution follows easily from the Lax-Milgram theorem by observing that $a_h(\cdot,\cdot)$ is continuous and coercive with respect to $\vert \cdot\vert_{2,h}$ , which is a norm in $H^{2,nc}$ in view of Lemma \ref{eq:norm2broken},  and thus on $V_{h,\ell}$ for any $\ell\geq 2$, cf. \eqref{global-space}.
We now address the proof of \eqref{vem-cea}. Using the triangular inequality we have that
$$| u-u_h |_{2,h} \leq | u-u^I |_{2,h} + | u_h - u^I|_{2,h}.$$
Setting $\delta_h=u_h-u^I$, employing \eqref{eq:stability} and \eqref{eq:1} we obtain the developments
\begin{eqnarray}
| \delta_h |_{2,h}^2 &=&  \sum_{K \in \Tau_h} a^K (\delta_h,\delta_h) 
\lesssim  \sum_{K \in \Tau_h} a_h^K (\delta_h,\delta_h)  
\nonumber\\
&=& a_h (u_h,\delta_h) - a_h (u^I,\delta_h) = 
\langle f_h, \delta_h \rangle - a_h (u^I,\delta_h)
\nonumber\\[1.em]
&=& \langle f_h, \delta_h \rangle - \sum_{K \in \Tau_h} a_h^K (u^I-u_\pi,\delta_h) -  
\sum_{K \in \Tau_h} a_h^K (u_\pi,\delta_h)
\nonumber\\
&=& \langle f_h, \delta_h \rangle - \sum_{K \in \Tau_h} a_h^K (u^I-u_\pi,\delta_h) -  
\sum_{K \in \Tau_h} a^K (u_\pi,\delta_h)
\nonumber\\
&=& \langle f_h, \delta_h \rangle - \sum_{K \in \Tau_h} a_h^K (u^I-u_\pi,\delta_h) + 
\sum_{K \in \Tau_h} a^K (u- u_\pi,\delta_h) - \sum_{K \in \Tau_h} a^K (u,\delta_h)
\nonumber\\
&=&\langle f_h-f, \delta_h \rangle + \N(u,\delta_h) - \sum_{K\in \Tau_h} a_h^K(u^I-u_\pi,\delta_h) + \sum_{K\in \Tau_h} a^K(u-u_\pi,\delta_h), \nonumber
\end{eqnarray}  
from which inequality~\eqref{def:N} follows.
\end{proof}
Finally, from the above result and bounding each term in \eqref{vem-cea} we obtain an estimate of the approximation error in the broken energy norm as stated in the following theorem.  
\begin{theorem}\label{th:conv}
Let us assume that the solution to \eqref{pb:w} satifies $u\in H^3(\Omega)$. Under the regularity assumption $({\mathbf H})$ on the mesh $\Tau_h$,  
for $\ell\geq 2$ the unique solution $u_h\in V_{h,\ell}$ to \eqref{pb:ncVEM} satisfies the following error estimate
\begin{equation}\label{error}
| u-u_h |_{2,h}\lesssim h^{\ell-1}.
\end{equation}
\end{theorem}
\begin{proof}
In order to prove \eqref{error} it is sufficient to combine \eqref{rhs:error} with \eqref{vem-cea}, use standard interpolation error estimates and the estimate of the conformity error $\N(u,\delta_h)$.
Let us focus on the last step. Assuming that $u$ is sufficiently smooth we rewrite the conformity error as follows:
\begin{align}
\N(u,\delta_h)
&\!=\! D\sum_{e\in \E_h}\left\{ \int_{e} (\Delta u - (1-\nu) u_{,tt}) {[\delta_h}_{,n}] ds  - \int_{e} \partial_n(\Delta u) [\delta_h] ds+ \int_e (1-\nu) u_{,nt} [{\delta_h}_{,t}] ds \right\}\nonumber\\
&\!=\! D(I+II+III).\nonumber
\end{align}
To estimate the above terms we employ the fact that $\delta_h$ belongs to $V_{h,\ell}$ and use standard interpolation error estimates for the $L^2$-projection $\Pi_e^{\ell}$ on polynomials defined on $e$. In particular, for the first term for $\ell\geq 2$ we use the definition of the global virtual element space $V_{h,\ell}$ and we find that
\begin{eqnarray}
I&=&\sum_{e\in \E_h} \int_{e} (I-\Pi_e^{\ell-2})(\Delta u - (1-\nu) u_{,tt}) {(I-\Pi_e^0)[\delta_h}_{,n}] ds \nonumber\\
&\lesssim& h^{\ell-2+1 -\frac 1 2}h^{\frac 1 2}   | \delta_h|_{2,h}= h^{\ell-1} | \delta_h |_{2,h}.\nonumber
\end{eqnarray}
For the second term we first consider the case $\ell=2,3$ and, \paa{in view of (H3)}, introduce, \paa{for each edge $e\subset \partial K$}, the linear Lagrange interpolant $I^1_{T^{(e)}}$ \paa{of $\delta _h$} on the triangle $T^{(e)}$, \paa{ which is obtained by connecting the point $x_B$ (interior to $K$) and the endpoints of $e$}. 
\paa{Clearly,  due to the $H^2$ regularity  of $\delta_h$, the interpolant $I^1_{T^{(e)}} \delta_h$ can be built based on employing the values of $\delta_h$ at the vertices of $T^{(e)}$}. In particular, using the fact that 
$\delta_h$ is continuous \paa{at the endpoints of $e$}  we have $[I^1_{T^{(e)}} \delta_h]_{\vert e}=0$. Hence, for $\ell=2$ employing standard interpolation error estimates and a trace inequality we get
\begin{eqnarray}
II &=& \sum_e \int_e \partial_n(\Delta u) ([\delta_h] - [I^1_{T^{(e)}} \delta_h]) ds \nonumber\\
&\lesssim&  h^{2-\frac 1 2}  | \delta_h|_{2,h}.\label{H3-eq-1}
\end{eqnarray}
On the other hand, for $\ell=3$, we have
\begin{eqnarray}
II &=& \sum_e \int_e (I-\Pi_e^0)(\partial_n(\Delta u)) ([\delta_h] - [I^1_{T^{(e)}} \delta_h]) ds \nonumber\\
&\lesssim&  h^{2}  | \delta_h|_{2,h},\label{H3-eq-2}
\end{eqnarray}
where we employed the definition of the global space $V_{h,3}$ together with standard interpolation error estimates and a trace inequality. In case $\ell \geq 4$ we have 
\begin{eqnarray}
II &=& \sum_e \int_e (I-\Pi_e^{\ell-3})(\partial_n(\Delta u)) (I-\Pi_e^1)[\delta_h]  ds \nonumber\\
&\lesssim& h^{\ell-3+1 -\frac 1 2}h^{1 +1 - \frac 1 2}   | \delta_h|_{2,h}=  h^{\ell-1}  | \delta_h|_{2,h}.
\end{eqnarray}
Finally, we consider the third term. Using the fact that $\delta_h$ is continuous at the vertexes \paa{of $\Tau_h$} and the fact that $\int_e [\delta_h] p ds =0 ~ \forall p\in \mathbb{P}^{\ell-3}(e)$ we deduce after integration by parts that  $\int_e  [{\delta_h}_{,t}]  q ds =0  ~ \forall q\in \mathbb{P}^{\ell-2}(e)$. Indeed, after observing that for any 
$p\in \mathbb{P}^{\ell-3}(e)$ there exists $q \in \mathbb{P}^{\ell-2}(e)\setminus \mathbb{P}^{0}(e)$
such that $p=q'$ we have 
\begin{eqnarray*}
0&=&\int_e [\delta_h] p ds =\int_e [\delta_h] q^\prime ds= - \int_e [\delta_{h,t}] q ds + ([\delta_{h}]q)(\vrtx_2)-([\delta_{h}]q)(\vrtx_1)\nonumber\\
&=&- \int_e [\delta_{h,t}] q,
\end{eqnarray*}
where we used the fact that \paa{the jump $[\delta_{h}]$ is zero when evaluated at the endpoints $\vrtx_1$ and $\vrtx_2$ of $e$}. 
Finally, for 
$q \in\mathbb{P}^{0}(e)$ we immediately have, after integration by parts,  
$$  \int_e [\delta_{h,t}] q =0. $$
In view of the above result we get 
\begin{eqnarray}
III &=& \sum_e \int_e (1-\nu) (I-\Pi_e^{\ell-2})u_{,nt} (I-\Pi_e^0)[{\delta_h}_{,t}] ds \nonumber\\
&\lesssim& h^{\ell-2+1 -\frac 1 2}h^{\frac 1 2}   | \delta_h|_{2,h}= h^{\ell-1} | \delta_h |_{2,h}\nonumber,
\end{eqnarray}
and this concludes the proof.
\end{proof}
%

\newcommand{\ilev}{n}
\newcommand{\nR}   {N_{P}}
\newcommand{\nF}   {N_{F}}
\newcommand{\nV}   {N_{V}}
\newcommand{\ndof} [1]{|V^{#1}_{h}|}
\newcommand{\hmax} {h}

\newcommand{\MeshONE}  {$\big\{\Tau_h^{(1)}\big\}_h$}
\newcommand{\MeshTWO}  {$\big\{\Tau_h^{(2)}\big\}_h$}
\newcommand{\MeshTHREE}{$\big\{\Tau_h^{(3)}\big\}_h$}
\newcommand{\MeshFOUR} {$\big\{\Tau_h^{(4)}\big\}_h$}


\section{Numerical results}\label{S:3}

The numerical experiments presented in this section are aimed to
confirm the a priori analysis developed in the previous
sections.
To study the accuracy of our new nonconforming method we solve 
the biharmonic problem \eqref{pb:1}-\eqref{pb:3} on the domain 
$\Omega=]0,1[\times]0,1[$.
The forcing term $f$ in~\eqref{pb:1} is set in accordance with the exact
solution:
\begin{align*}
  u(x,y) = x^2(1-x)^2\,y^2(1-y)^2,
\end{align*}
which obviously satisfies the boundary conditions in~\eqref{pb:2}-\eqref{pb:3}.\\

The performance of the VEM is investigated by observing experimentally
the convergence behavior on four different sequences of unstructured meshes
\paa{labelled} by~\MeshONE{}, \MeshTWO{}, \MeshTHREE{}, and \MeshFOUR{}.
All mesh data are reported in
Tables~\ref{table:mesh-data:criss-cross}-\ref{table:mesh-data:randomized:quadrilateral} \paa{in the Appendix}.
Figs.~\ref{fig:meshes}(a)-(d) show the first and
second mesh of each sequence (top and right plots, respectively).
The meshes in~\MeshONE{}, also known as \emph{criss-cross meshes}, are composed by
first partitioning $\Omega$ in regular square grids and then splitting 
each square cell into four triangular subcells by connecting the four vertices 
along the diagonal.
It is worth \paa{recalling} that our nonconforming VEM for $\ell=2$ 
on triangular meshes coincides with the Morley finite element 
method~\cite{Morley:1968}. 
The meshes in~\MeshTWO{} are built as follows.
First, we determine a primal mesh by remapping the position
$(\hat{x},\hat{y})$ of the nodes of a uniform square partition of
$\Omega$ by the smooth coordinate transformation (see
~\cite{Brezzi-Lipnikov-Shashkov:2005}):
\begin{align*}
  x &= \hat{x} + 0.1\, \sin(2\pi\hat{x})\sin(2\pi\hat{y}),\\
  y &= \hat{y} + 0.1\, \sin(2\pi\hat{x})\sin(2\pi\hat{y}).
\end{align*}
Then, the corresponding mesh of \MeshTWO{} is built from the primal mesh by
splitting each quadrilateral cell into two triangles and connecting
the barycenters of adjacent triangular cells by a straight segment.
The mesh construction is completed at the boundary by connecting the
barycenters of the triangular cells close to the boundary to the
midpoints of the boundary edges and these latter ones to the boundary
vertices of the primal mesh.
The meshes in \MeshTHREE{} are obtained by filling the unit square
with a suitably scaled non-convex octagonal cell, which is cut at the
domain boundaries to fit into the unit square domain.
The meshes in \MeshFOUR{} are built by partitioning the domain $\Omega$
into square cells and relocating each interior node to a random
position inside a square box centered at that node.
The sides of this square box are aligned with the coordinate axis and
their length is equal to $0.8$ times the minimum distance between two
adjacent nodes of the initial square mesh.
All the meshes are parametrized by the number of partitions in each
direction.
The starting mesh of every sequence is built from a $5\times 5$
regular grid, while for the $\ilev$-th refined mesh the underlying resolution 
is $10\ilev\times 10\ilev$.
For the virtual element spaces of order $2\leq\ell\leq 4$ we consider a sequence 
of 9 meshes; for $\ell=5$ the calculation is arrested after the fifth mesh when
rounding errors begin affecting the accuracy of the approximation due to the
increasing ill-conditioning of the \paa{algebraic} problem.\\

\begin{figure}[thb]
  \centering
  \begin{tabular}{cc}
    \begin{overpic}[width=0.8\textwidth]{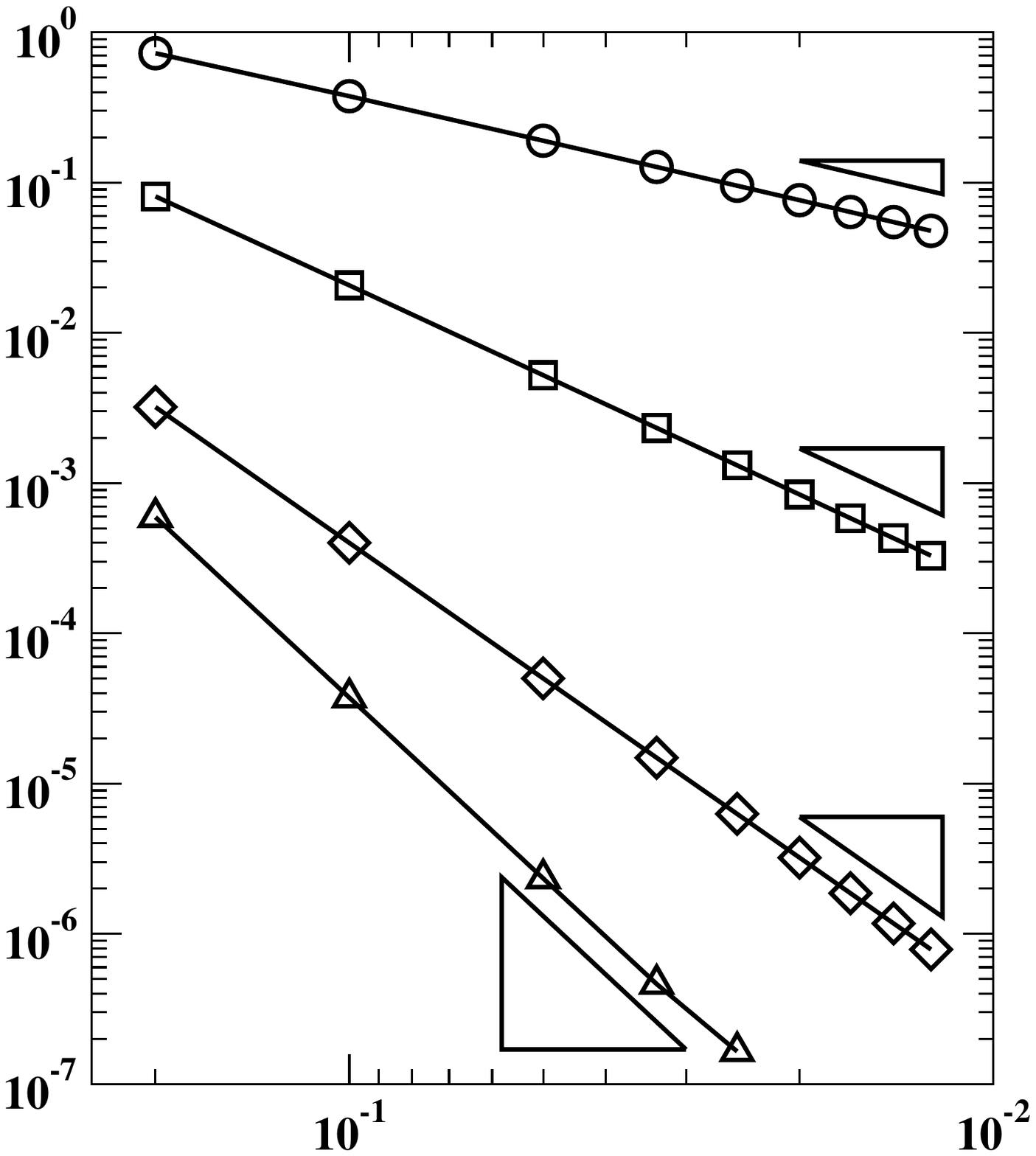}
    \put(27,2){\textbf{Mesh size $h$}}
    \put(2,17.5){\begin{sideways}\textbf{ 2h Approximation Error }\end{sideways}}
    \put(51,60.5){\textbf{1}}
    \put(51,45.5){\textbf{2}}
    \put(51,26.5){\textbf{3}}
    \put(30,16.5){\textbf{4}}
    \end{overpic}
    &\hspace{-5cm}
    \begin{overpic}[width=0.8\textwidth]{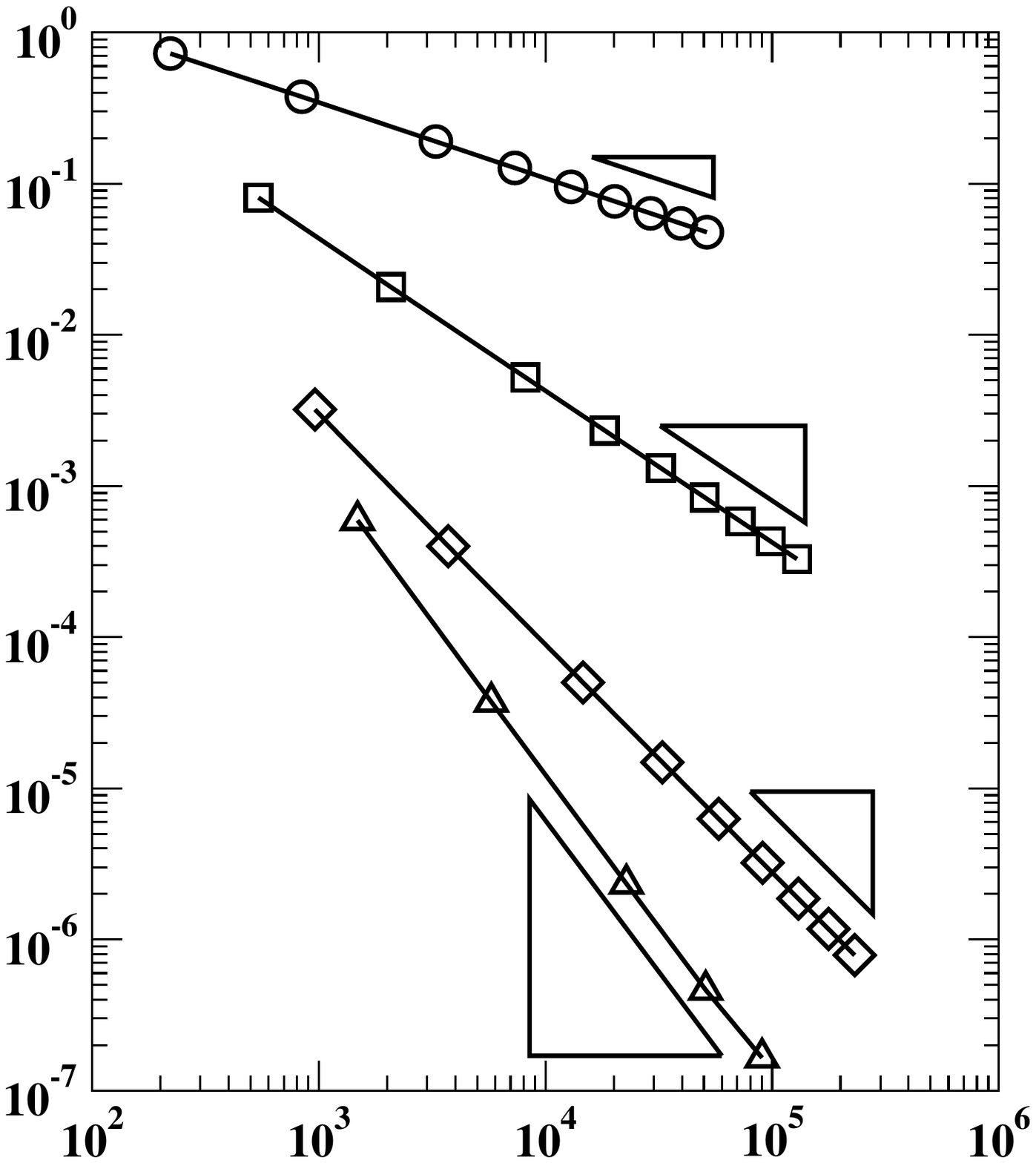}
    \put(30,2){\textbf{\#dofs}}
    \put(2,17.5){\begin{sideways}\textbf{ 2h Approximation Error }\end{sideways}}
    \put(44,57.5){$\mathbf{\frac{1}{2}}$}
    \put(49,42.5){$\mathbf{1}$}
    \put(52,22.5){$\mathbf{\frac{3}{2}}$}
    \put(31,18.5){$\mathbf{2}$}
    \end{overpic}
  \end{tabular}
  \caption{Relative $2h$-approximation errors using the sequence of criss-cross triangular meshes
  versus the mesh size parameter $h$ (left panel) and the total number of degrees of freedom 
  $\#$dofs (right panel).}
  \label{fig:errors:M220}


  \centering
  \begin{tabular}{cc}
    \begin{overpic}[width=0.8\textwidth]{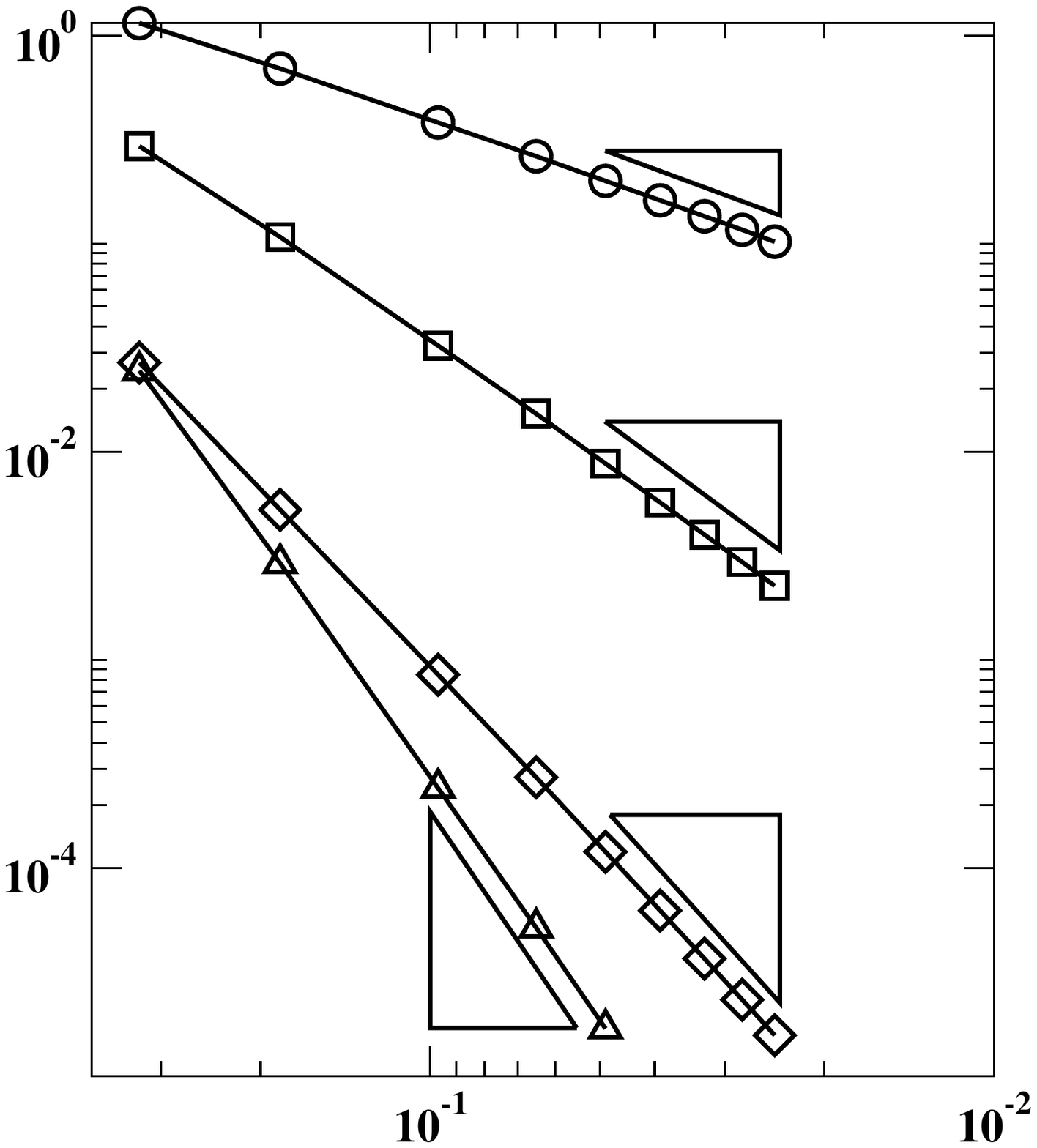}
    \put(27,2){\textbf{Mesh size $h$}}
    \put(2,17.5){\begin{sideways}\textbf{ 2h Approximation Error }\end{sideways}}
    \put(48,56.5){\textbf{1}}
    \put(48,40.5){\textbf{2}}
    \put(48,18.5){\textbf{3}}
    \put(26,18.5){\textbf{4}}
    \end{overpic}
    &\hspace{-5cm}
    \begin{overpic}[width=0.8\textwidth]{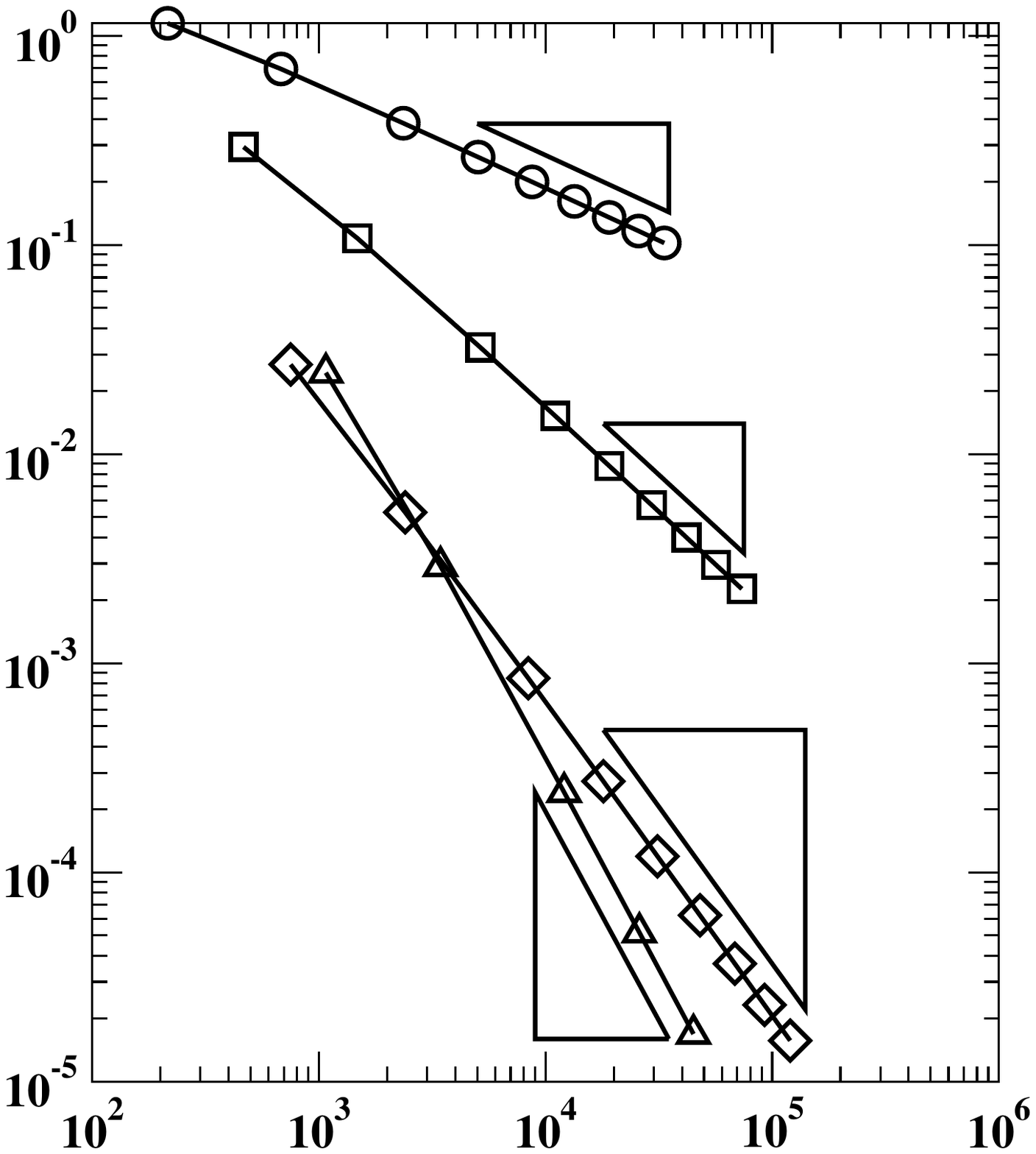}
    \put(30,2){\textbf{\#dofs}}
    \put(2,17.5){\begin{sideways}\textbf{ 2h Approximation Error }\end{sideways}}
    \put(42,  57.5){$\mathbf{\frac{1}{2}}$}
    \put(45.5,40.5){$\mathbf{1}$}
    \put(49,  22.5){$\mathbf{\frac{3}{2}}$}
    \put(31,  18.5){$\mathbf{2}$}
    \end{overpic}
    \\[-0.5em]
  \end{tabular}
  \caption{Relative $2h$-approximation errors using the sequence of remapped hexagonal meshes
  versus the mesh size parameter $h$ (left panel) and the total number of degrees of freedom 
  $\#$dofs (right panel).}
  \label{fig:errors:Md201}
\end{figure}

\begin{figure}[t]
  \centering
  \begin{tabular}{cc}
    \begin{overpic}[width=0.8\textwidth]{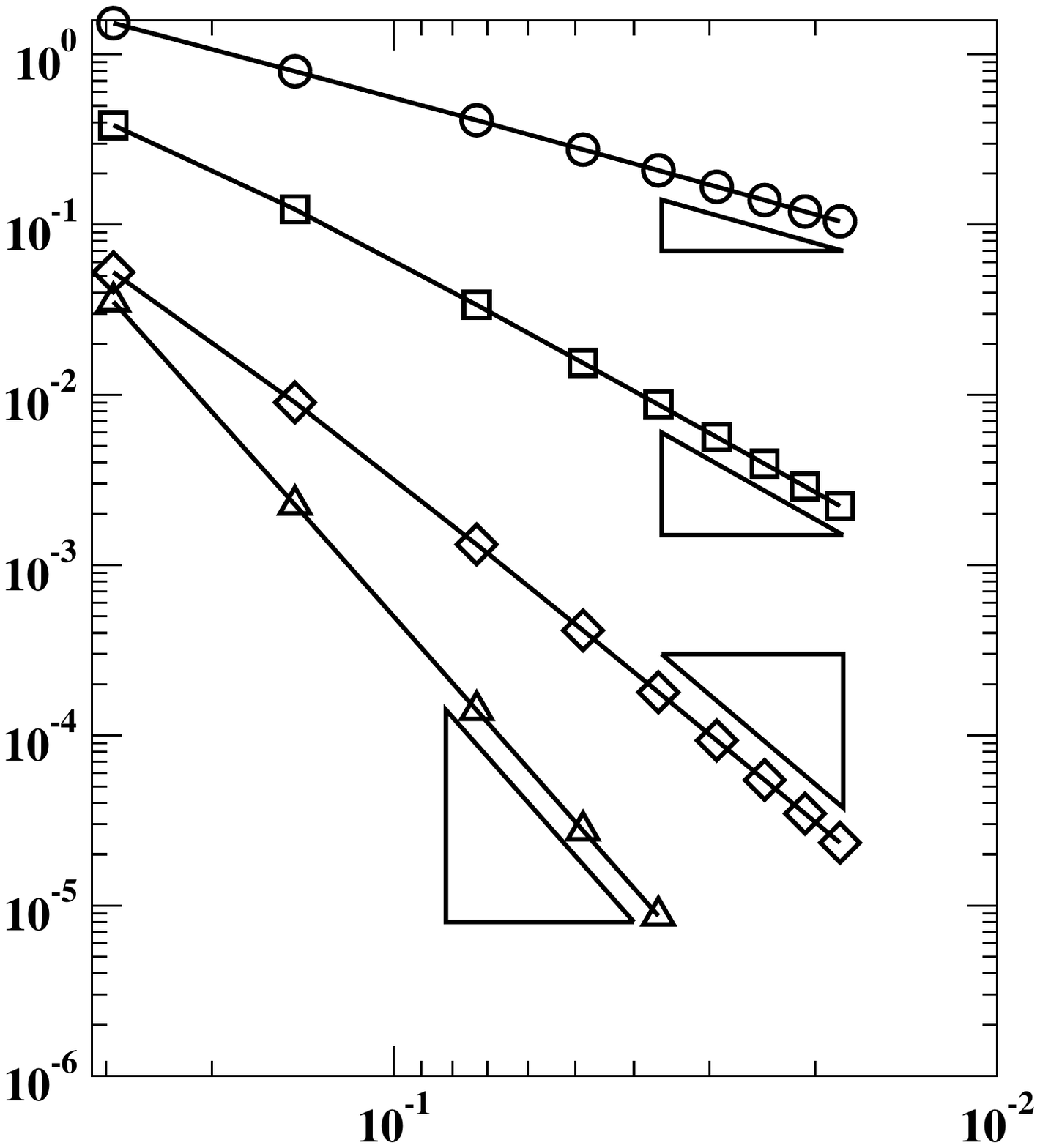}
    \put(27,2){\textbf{Mesh size $h$}}
    \put(2,17.5){\begin{sideways}\textbf{ 2h Approximation Error }\end{sideways}}
    \put(38,54.5){\textbf{1}}
    \put(38,41.0){\textbf{2}}
    \put(51,28.5){\textbf{3}}
    \put(27,23.5){\textbf{4}}
    \end{overpic}
    &\hspace{-5cm}
    \begin{overpic}[width=0.8\textwidth]{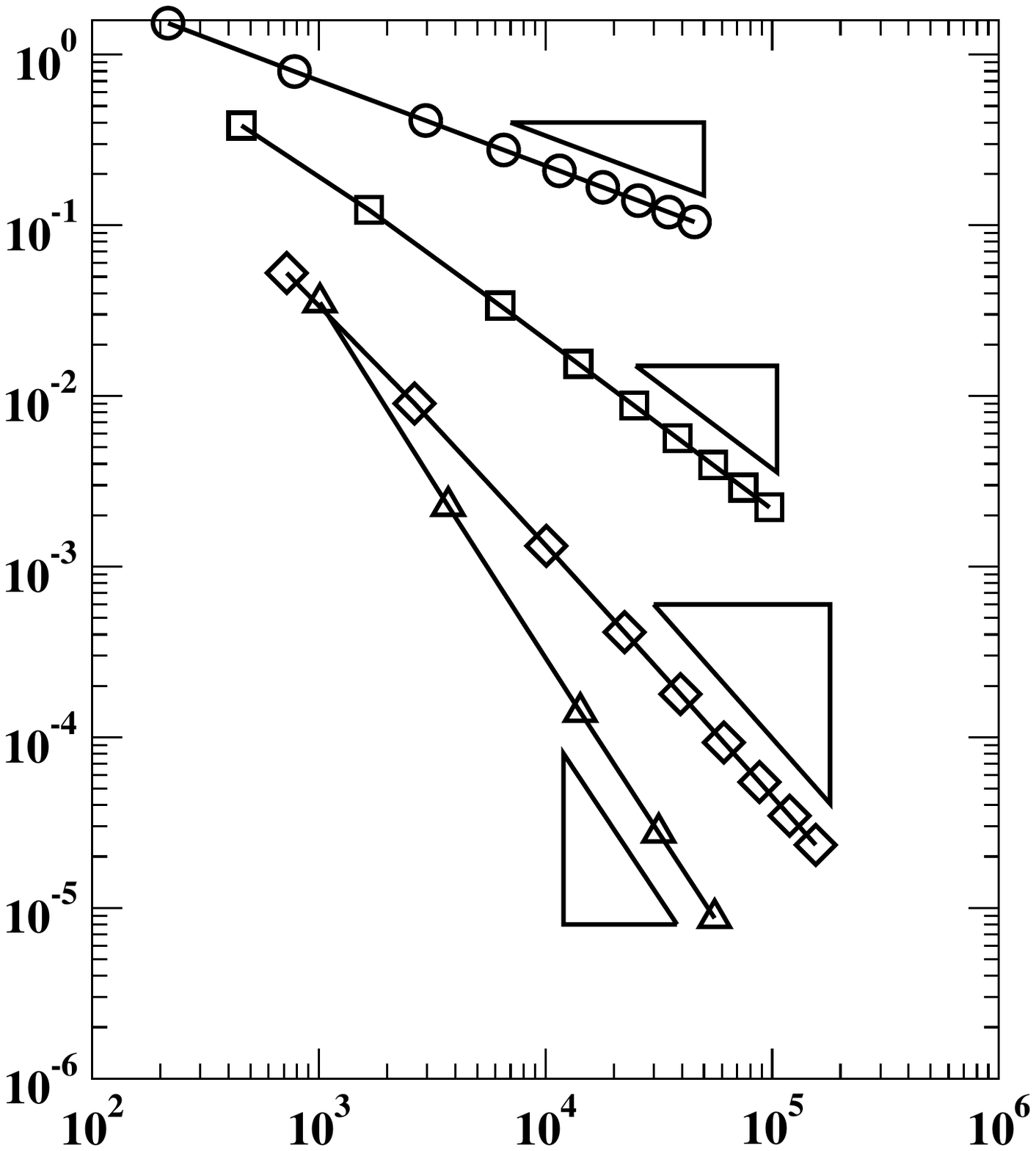}
    \put(30,2){\textbf{\#dofs}}
    \put(2,17.5){\begin{sideways}\textbf{ 2h Approximation Error }\end{sideways}}
    \put(43,  57.5){$\mathbf{\frac{1}{2}}$}
    \put(47.5,44.0){$\mathbf{1}$}
    \put(50,  29.5){$\mathbf{\frac{3}{2}}$}
    \put(32.5,22.5){$\mathbf{2}$}
    \end{overpic}
    \\[-0.5em]
  \end{tabular}
  \caption{Relative $2h$-approximation errors using the sequence of non-convex regular meshes
  versus the mesh size parameter $h$ (left panel) and the total number of degrees of freedom 
  $\#$dofs (right panel).}
  \label{fig:errors:M901}


  \centering
  \begin{tabular}{cc}
    \begin{overpic}[width=0.8\textwidth]{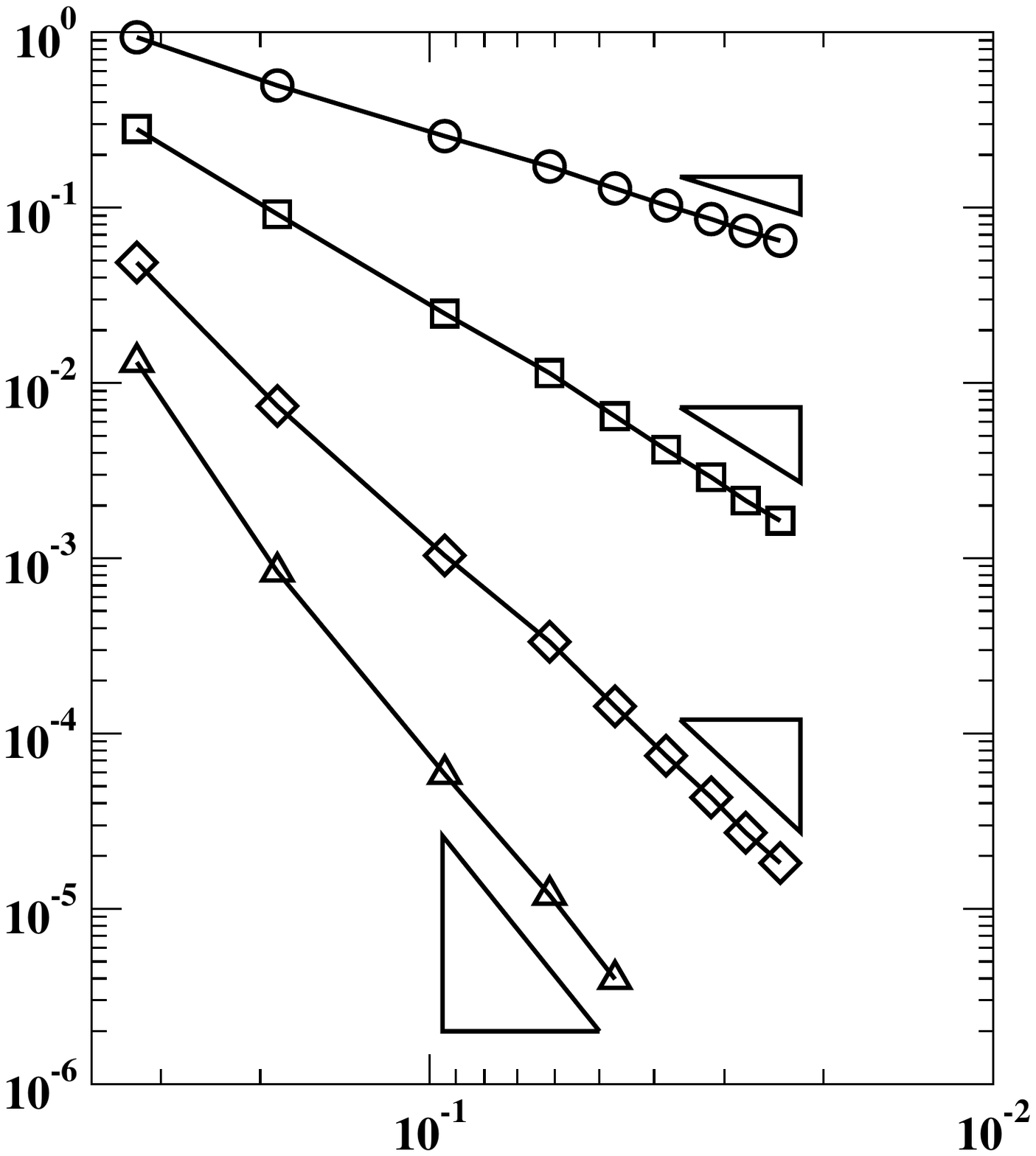}
    \put(27,2){\textbf{Mesh size $h$}}
    \put(2,17.5){\begin{sideways}\textbf{ 2h Approximation Error }\end{sideways}}
    \put(49,  56.5){\textbf{1}}
    \put(49,  43.5){\textbf{2}}
    \put(49,  26.5){\textbf{3}}
    \put(26.5,17.5){\textbf{4}}
    \end{overpic}
    &\hspace{-5cm}
    \begin{overpic}[width=0.8\textwidth]{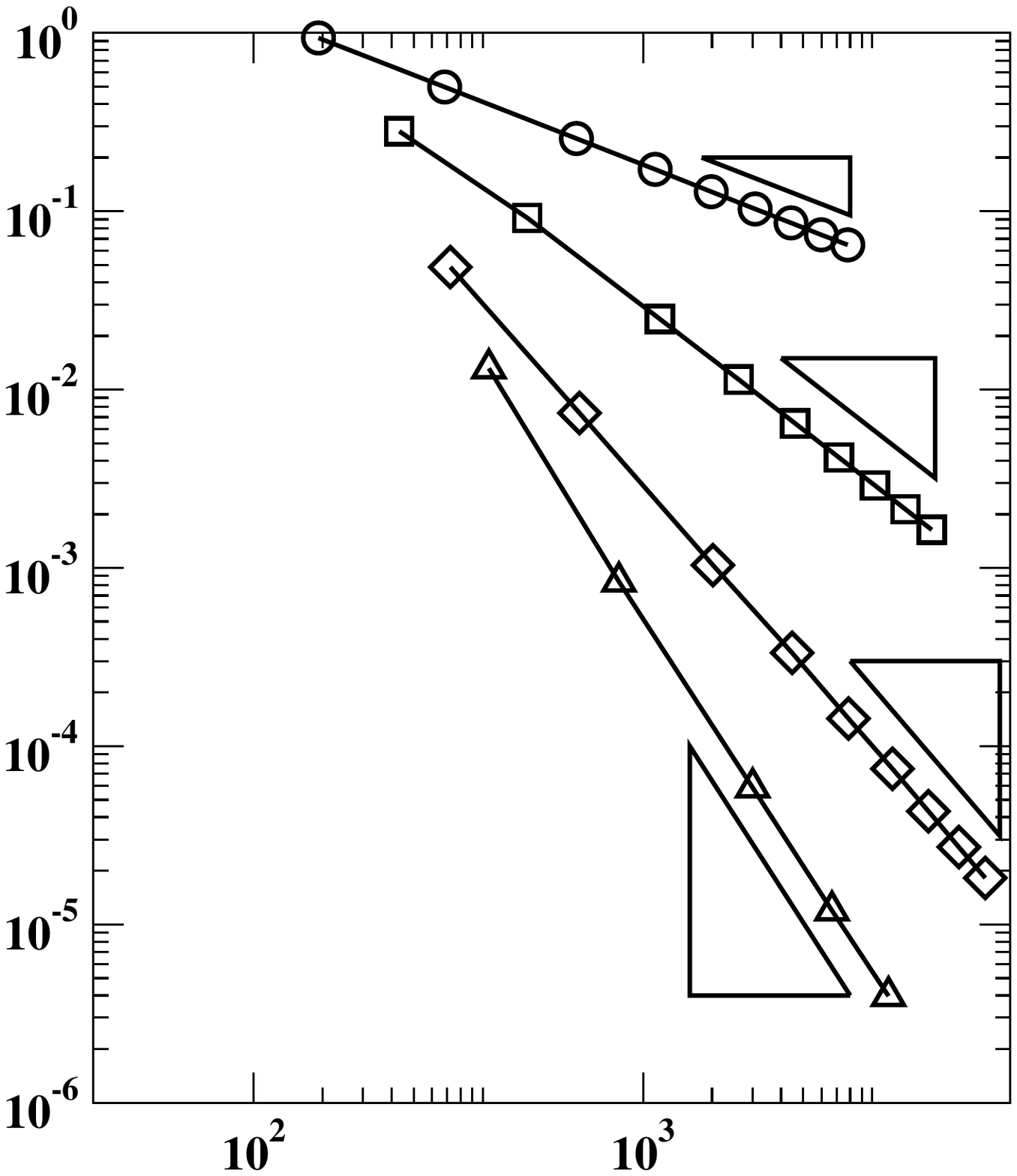}
    \put(30,2){\textbf{\#dofs}}
    \put(2,17.5){\begin{sideways}\textbf{ 2h Approximation Error }\end{sideways}}
    \put(50.5,57.0){$\mathbf{\frac{1}{2}}$}
    \put(54.5,45.5){$\mathbf{1}$}
    \put(53,  30.5){$\mathbf{\frac{3}{2}}$}
    \put(38.5,22.0){$\mathbf{2}$}
    \end{overpic}
    \\[-0.5em]
  \end{tabular}
  \caption{Relative $2h$-approximation errors using the sequence of randomized quadrilateral meshes
  versus the mesh size parameter $h$ (left panel) and the total number of degrees of freedom 
  $\#$dofs (right panel).}
  \label{fig:errors:M102}
\end{figure}

For $\ell\geq 2$, we define the relative ``$2h$'' error by
\begin{align*}
\textsf{Error}_{2,h}=\paa{\frac{ | \Pi^{\Delta}_{\ell}(u-u_h) |_{2,h} }{ | \Pi^{\Delta}_\ell(u) |_{2,h} }.}
\end{align*}
\paa{with $ {\Pi^{\Delta}_\ell}_{\vert K}= \Pi^{\Delta,K}_{\ell}$.}
Thus, on every element $K\in\Tau_h$, we compare \paa{$\Pi^{\Delta,K}_\ell u$, 
the elliptic projection
of the
exact
solution $u$ and
$\Pi^{\Delta,K}_{\ell}u_h$,} the projection of the virtual element solution
$u_h$.
These relative errors are shown in the log-log plots of 
Figs.~\ref{fig:errors:M220},
\ref{fig:errors:Md201},
\ref{fig:errors:M901}, and
\ref{fig:errors:M102}, with respect to the mesh size parameter
$h$ (left panels) and the number
of degrees of freedom (right panels).
The convergence rate is reflected by the slope of the experimental
error curves that are obtained by joining the error values
measured on the sequence of refined meshes for each polynomial degree
$2\leq\ell\leq 5$.
Each experimental slope has to be compared with the theoretical slope, which is
shown for each curve by a triangle and whose value is indicated by the
nearby number.
From the a priori analysis of Section~\ref{S:error}, cf.
Theorem~\ref{th:conv} and inequality~\eqref{error}), the $2h$-approximation errors must decrease proportionally to $h^{\ell-1}$
when we use the virtual element space $V^{h}_{\ell}$.
These errors are also expected to decrease proportionally to 
$\ndof{\ell}^{-\frac{\ell-1}{2}}$, where $\ndof{\ell}$ is the total number of
degrees of freedom of the $\ell$-th virtual element space, 
because $\ndof{\ell}$ is roughly proportional to $h^{-{1}\slash{2}}$,
Accordingly, the experimental slopes for $\textsf{Error}_{2,h}$
are expected to be closed to $\ell-1$ and $(\ell-1)\slash{2}$ when 
we plot the error curves versus the mesh size parameter $h$ and the number of
degrees of freedom.
The experimental convergence rates are in perfect agreement with the
theoretical ones for all such calculations.\\

Finally, it is worth mentioning that in a preliminary stage of this
work, the consistency of the nonconforming VEM of order $\ell$ for
$2\leq\ell\leq 5$, i.e., the
exactness of the method for polynomial solutions of degree up to $\ell$, has
been tested numerically by solving the bi-harmonic
equation~\eqref{pb:1} with Dirichlet boundary conditions and forcing term
determined by the monomials $x^{\mu}y^{\nu}$ for all possible combinations
of integers $\mu$ and $\nu$ such that $\mu+\nu\leq\ell$.
Non-homogeneous Dirichlet conditions were imposed in strong form by setting
the boundary degrees of freedom to the values determined by the exact solution.
For these experiments, we considered a wider set of polygonal meshes
(including the four considered in this section).
In all the cases, the magnitude of the $2h$ errors was comparable to the arithmetic precision, thus confirming the polynomial consistency of the method.
We also verified that our nonconforming VEM 
for $\ell=2$ on the criss-cross triangular meshes provides the same results of an 
independent implementation of the Morley finite element method.
For the sake of brevity, these results are not reported here.



\begin{figure}
\centering
\begin{tabular}{cccc}
\includegraphics[scale=0.135]{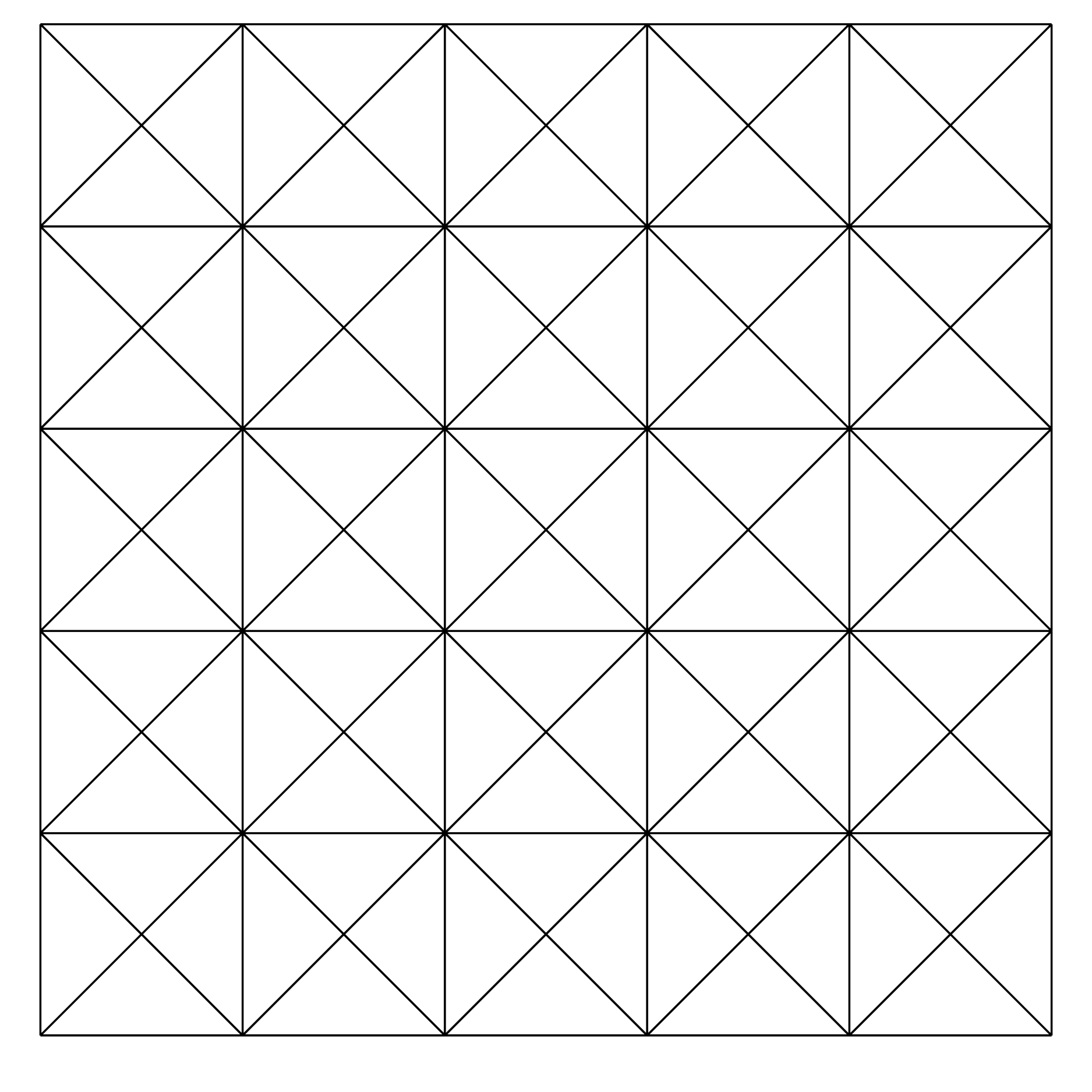}&
\includegraphics[scale=0.135]{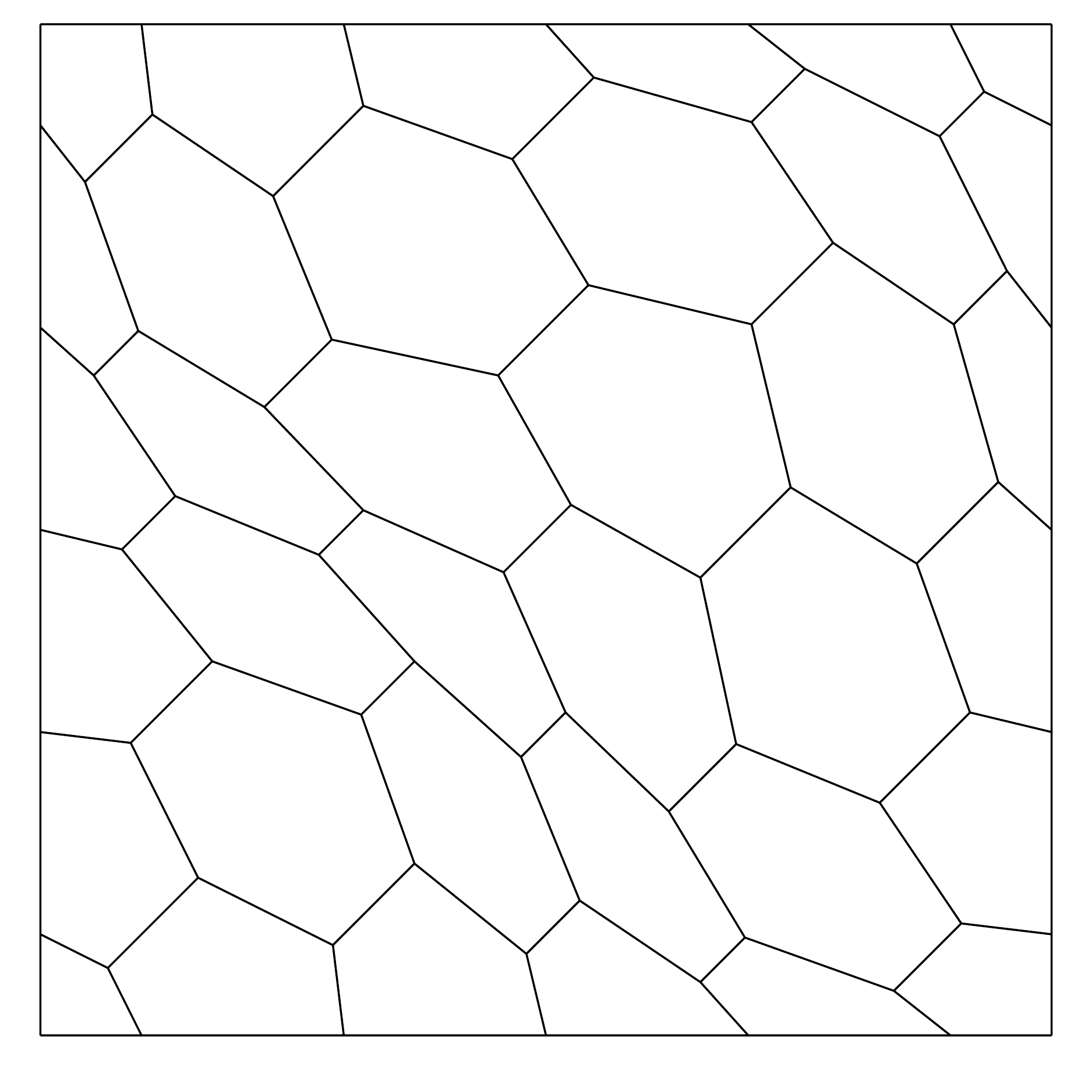}&
\includegraphics[scale=0.135]{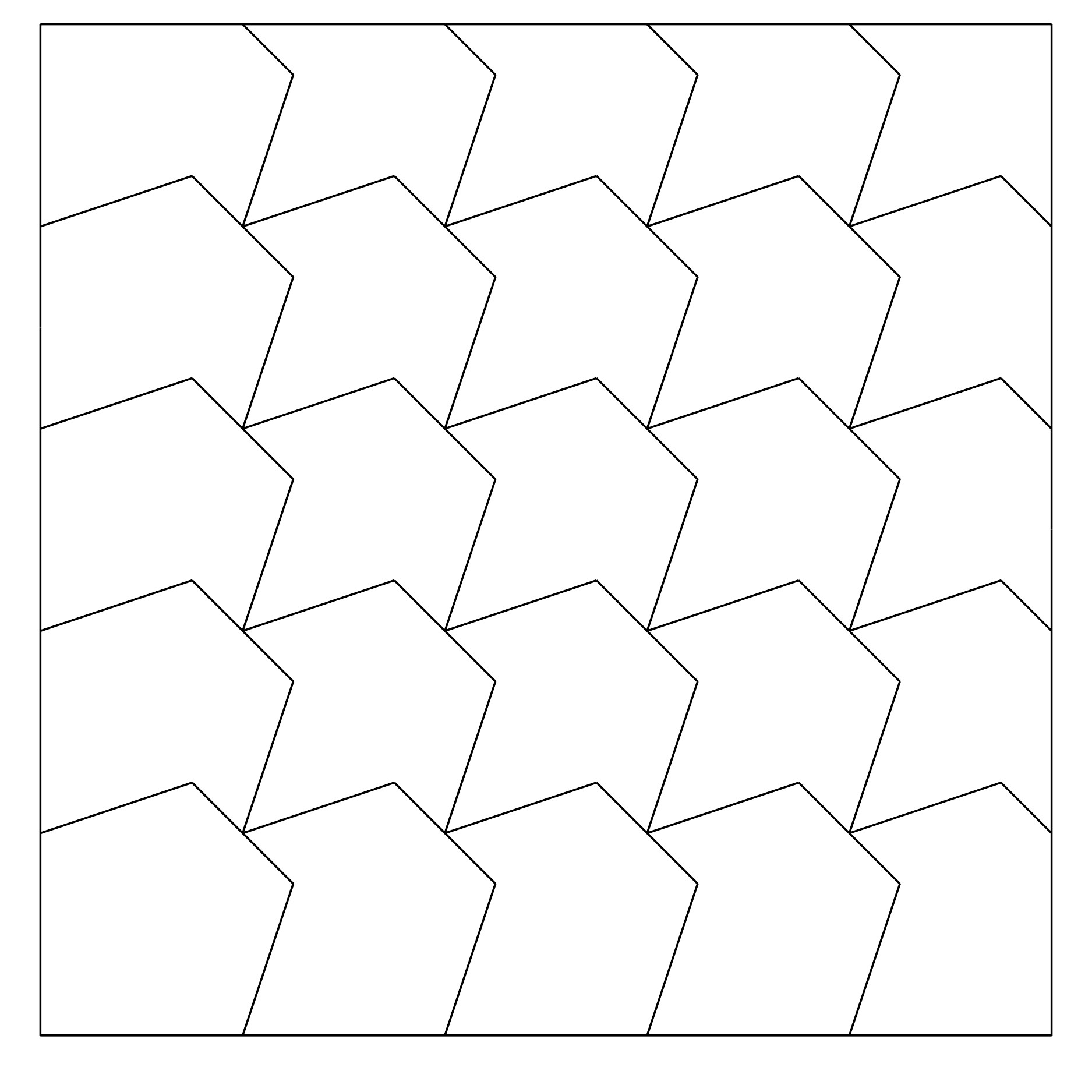}&
\includegraphics[scale=0.135]{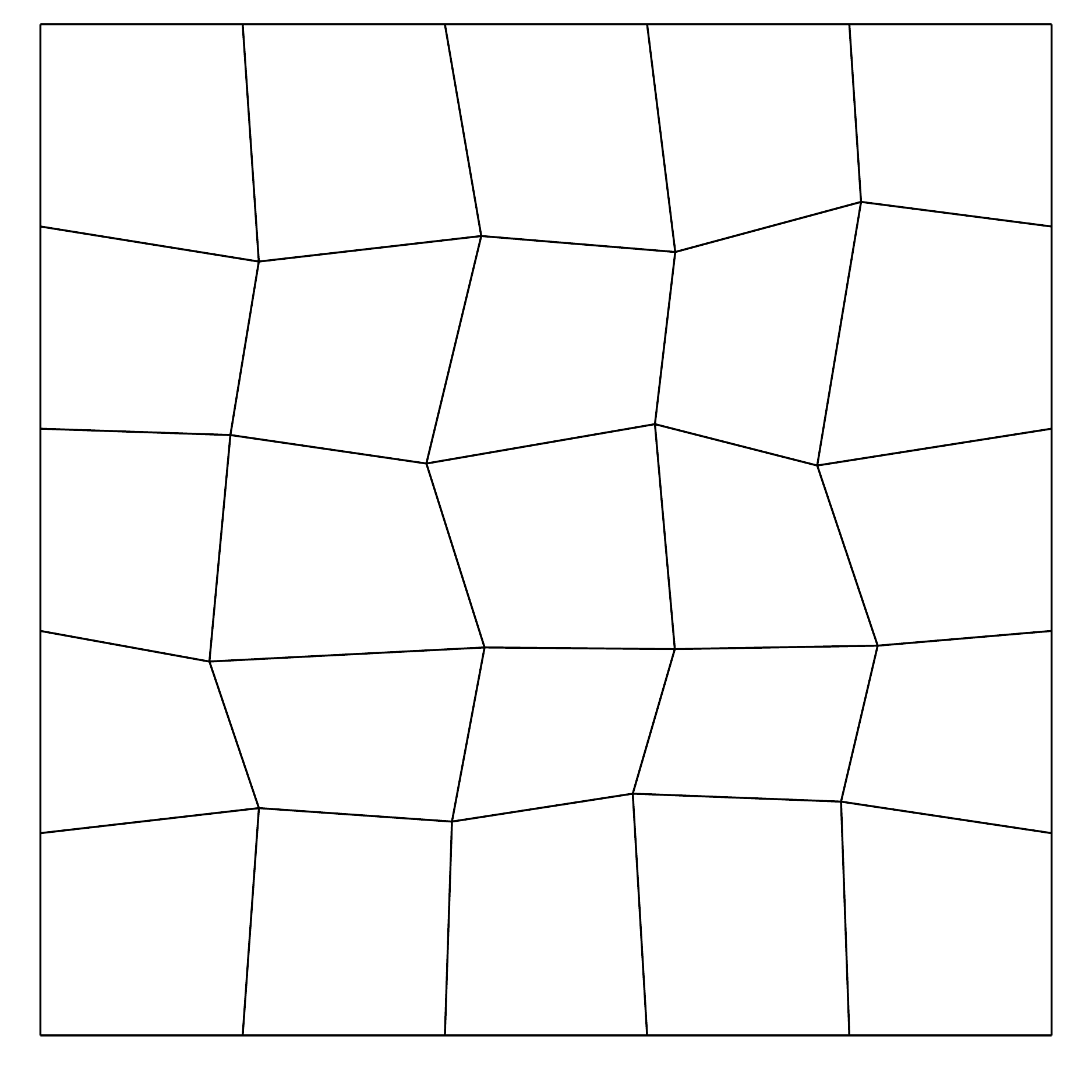}\\[1em]
\includegraphics[scale=0.135]{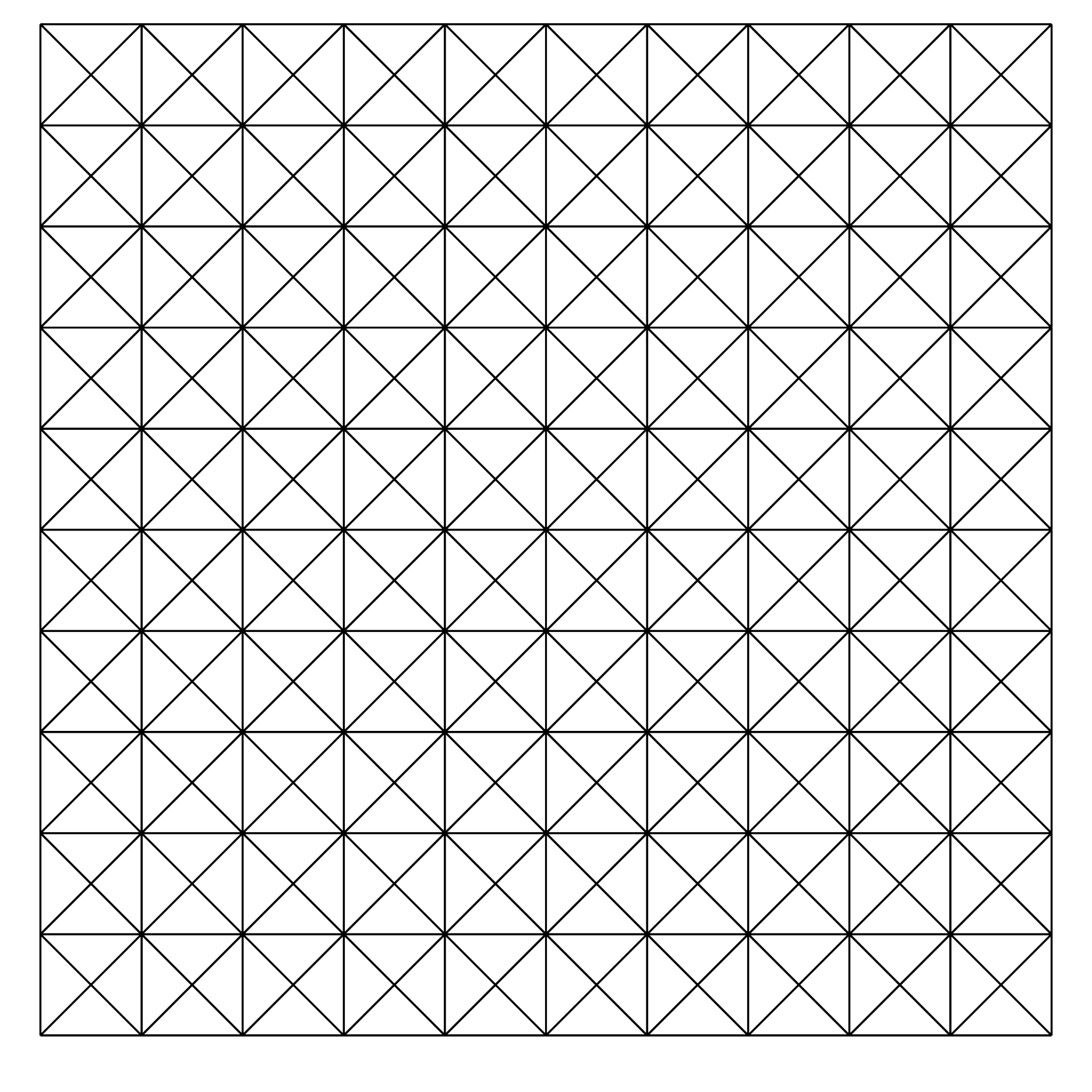}&
\includegraphics[scale=0.135]{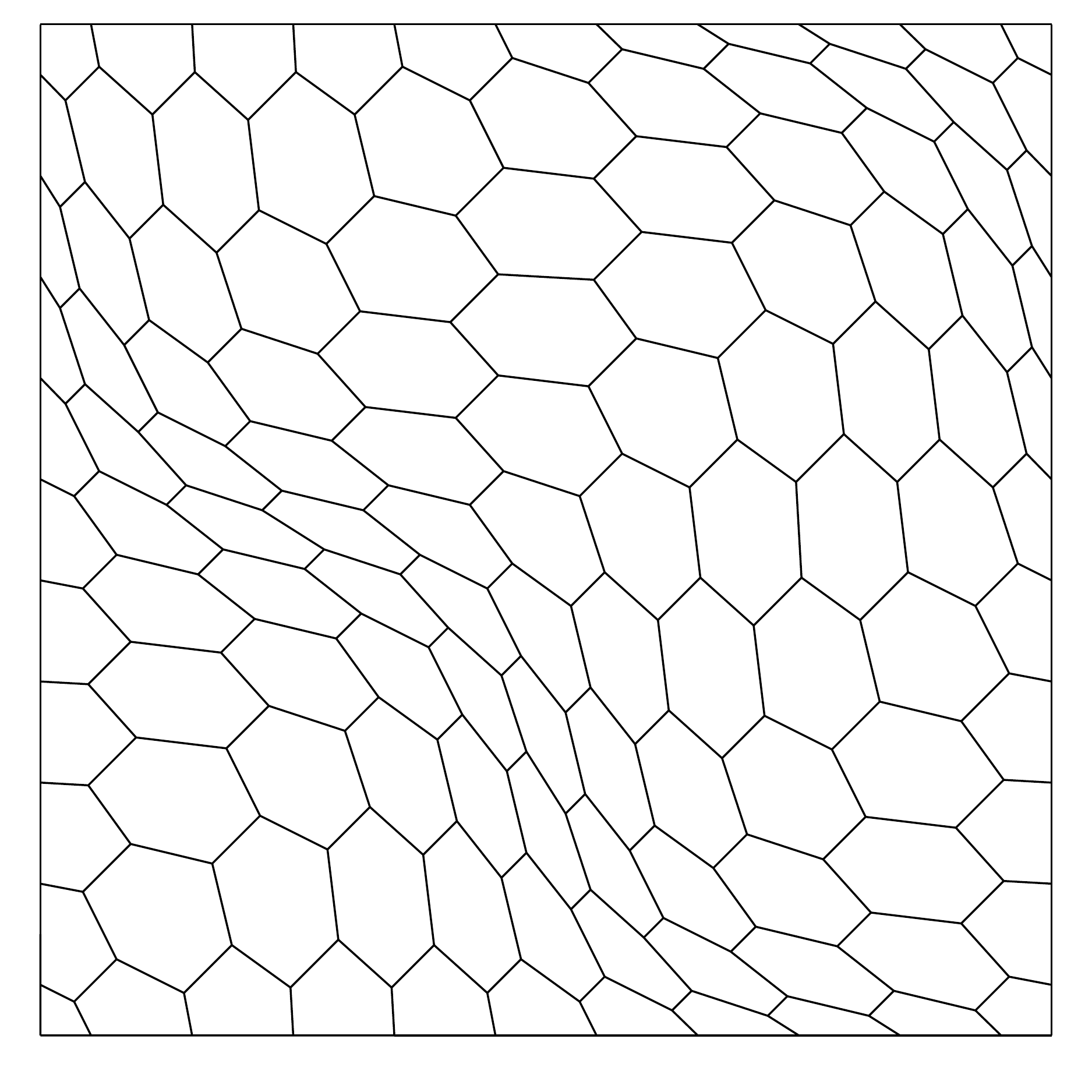}&
\includegraphics[scale=0.135]{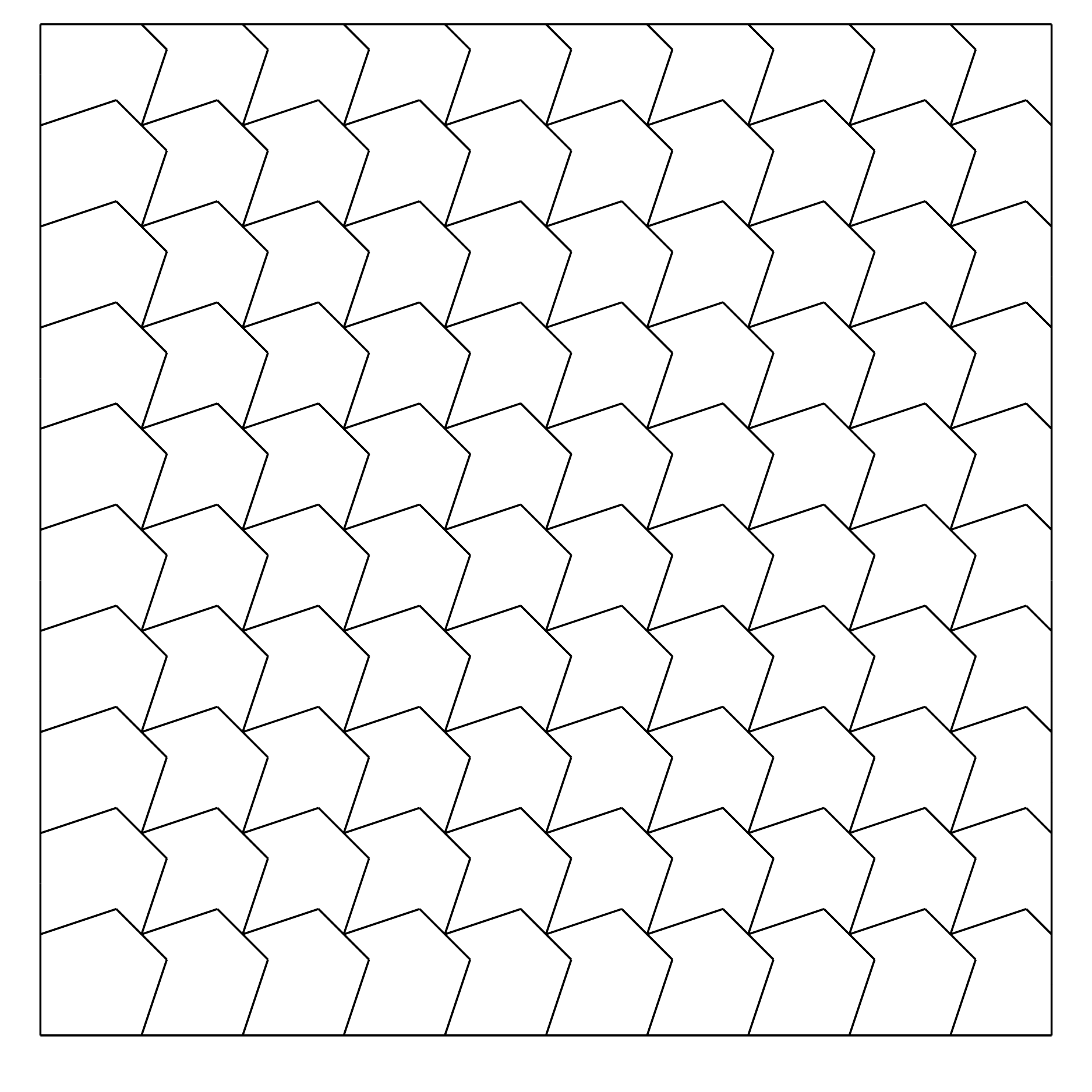}&
\includegraphics[scale=0.135]{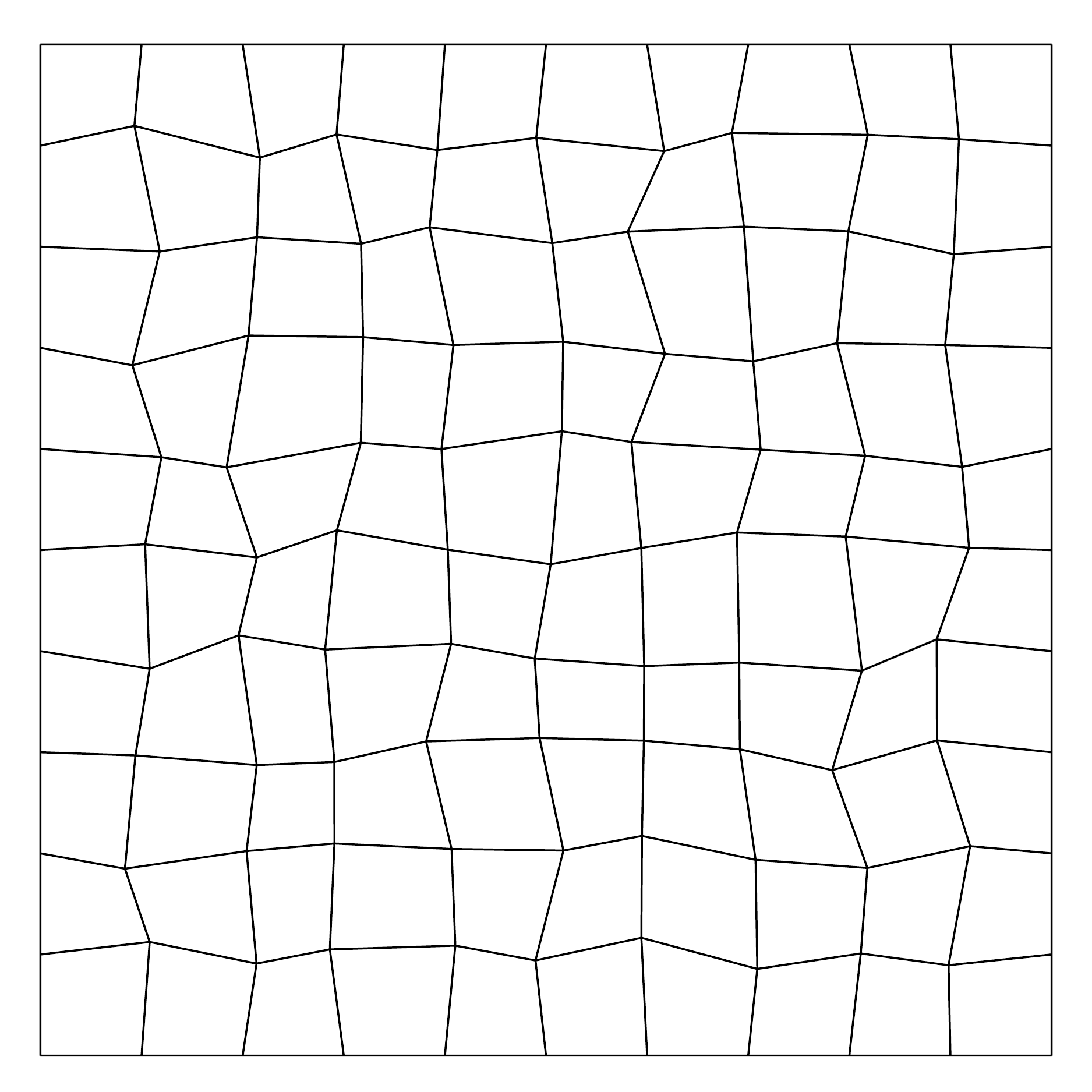}\\[1em]
\text{(a)} & \text{(b)} & \text{(c)} & \text{(d)}
\end{tabular}
\caption{Base mesh (top row) and first refinement (bottom row) of the four mesh families: (a) criss-cross triangular mesh; (b) mainly hexagonal mesh; (c) non-convex regular mesh; (d) randomized quadrilateral mesh.}
\label{fig:meshes}
\end{figure}


\section{Conclusions}\label{S:4}
In this paper we presented the arbitrary-order accurate fully nonconforming 
virtual element method for biharmonic problems on polygonal meshes.  
The virtual element space is made of functions that may be globally {\em not}-continuous. 
An optimal error estimate in the broken energy norm is derived for all polynomial 
approximation orders and numerical results assess the validity of the theoretical estimate.

\bibliographystyle{plain}
\bibliography{ncvem-biharm}

\clearpage

\section*{Appendix}
In Tables \ref{table:mesh-data:criss-cross}--\ref{table:mesh-data:randomized:quadrilateral},
we report the geometric data and the total number of degrees of freedom 
of the associated VEM spaces for the sequences of considered meshes.
More precisely, in Tables \ref{table:mesh-data:criss-cross}--\ref{table:mesh-data:randomized:quadrilateral}
the first column  reports the refinement level $\ilev=0,1,2,\ldots$,
the second, third and fourth columns show the corresponding total number of polygonal cells ($\nR$), faces ($\nF$) and vertexes ($\nV$), respectively, whereas in the fifth column the corresponding mesh size $\hmax$ is shown.
Finally, in the last four columns we report the total number of degrees of freedom of the corresponding VEM spaces
$V^\ell_h$, $\ell=2,\ldots,5$. 
\begin{table}[h]
\begin{small}
\caption{Geometric data and number of degrees of freedom of the sequence of criss-cross meshes.}
\label{table:mesh-data:criss-cross}
\begin{center}
\begin{tabular}{c|ccc|c|cccc}
$\ilev$ & $\nR$   & $\nF$   & $\nV$   & $\hmax$         & $\ndof{2}$ & $\ndof{3}$ & $\ndof{4}$ & $\ndof{5}$ \\
\hline
$0$     & $100$   & $160$   & $61$    & $2.00\,10^{-1}$ & $221$      & $541$      & $961$      & $1481$     \\
$1$     & $400$   & $620$   & $221$   & $1.00\,10^{-1}$ & $841$      & $2081$     & $3721$     & $5761$     \\
$2$     & $1600$  & $2440$  & $841$   & $5.00\,10^{-2}$ & $3281$     & $8161$     & $14641$    & $22721$    \\
$3$     & $3600$  & $5460$  & $1861$  & $3.33\,10^{-2}$ & $7321$     & $18241$    & $32761$    & $50881$    \\
$4$     & $6400$  & $9680$  & $3281$  & $2.50\,10^{-2}$ & $12961$    & $32321$    & $58081$    & $90241$    \\
$5$     & $10000$ & $15100$ & $5101$  & $2.00\,10^{-2}$ & $20201$    & $50401$    & $90601$    & $--$       \\ 
$6$     & $14400$ & $21720$ & $7321$  & $1.67\,10^{-2}$ & $29041$    & $72481$    & $130321$   & $--$       \\ 
$7$     & $19600$ & $29540$ & $9941$  & $1.43\,10^{-2}$ & $39481$    & $98561$    & $177241$   & $--$       \\ 
$8$     & $25600$ & $38560$ & $12961$ & $1.25\,10^{-2}$ & $51521$    & $128641$   & $231361$   & $--$       \\ 
\end{tabular}
\end{center}
\end{small}
\end{table} 
\begin{table}[h]
\begin{small}
\caption{Geometric data and number of degrees of freedom of the sequence of remapped hexagonal meshes.}
\label{table:mesh-data:remapped-hexagonal}
\begin{center}
\begin{tabular}{c|ccc|c|cccc}
$\ilev$ & $\nR$  & $\nF$   & $\nV$   & $\hmax$         & $\ndof{2}$ & $\ndof{3}$ & $\ndof{4}$ & $\ndof{5}$ \\
\hline
$0$     & $36$   & $125$   & $90$    & $3.28\,10^{-1}$ & $215$      & $465$      & $751$       & $1073$     \\
$1$     & $121$  & $400$   & $280$   & $1.85\,10^{-1}$ & $680$      & $1480$     & $2401$      & $3443$     \\
$2$     & $441$  & $1400$  & $960$   & $9.69\,10^{-2}$ & $2360$     & $5160$     & $8401$      & $12083$    \\
$3$     & $961$  & $3000$  & $2040$  & $6.49\,10^{-2}$ & $5040$     & $11040$    & $18001$     & $25923$    \\
$4$     & $1681$ & $5200$  & $3520$  & $4.89\,10^{-2}$ & $8720$     & $19120$    & $31201$     & $44963$    \\
$5$     & $2601$ & $8000$  & $5400$  & $3.91\,10^{-2}$ & $13400$    & $29400$    & $48001$     & $--$       \\ 
$6$     & $3721$ & $11400$ & $7680$  & $3.26\,10^{-2}$ & $19080$    & $41880$    & $68401$     & $--$       \\ 
$7$     & $5041$ & $15400$ & $10360$ & $2.80\,10^{-2}$ & $25760$    & $56560$    & $92401$     & $--$       \\ 
$8$     & $6561$ & $20000$ & $13440$ & $2.45\,10^{-2}$ & $33440$    & $73440$    & $120001$    & $--$       \\ 
\end{tabular}
\end{center}
\end{small}
\end{table}
\begin{table}[b]
\begin{small}
\caption{Geometric data and number of degrees of freedom of the sequence of nonconvex octagonal meshes.}
\label{table:mesh-data:nonconvex-octagonal}
\begin{center}
\begin{tabular}{c|ccc|c|cccc}
$\ilev$ & $\nR$  & $\nF$   & $\nV$   & $\hmax$         & $\ndof{2}$ & $\ndof{3}$ & $\ndof{4}$ & $\ndof{5}$ \\
\hline
$0$     & $25$   & $120$   & $96$    & $2.91\,10^{-1}$ & $216$     & $456$       & $721$      & $1011$     \\
$1$     & $100$  & $440$   & $341$   & $1.46\,10^{-1}$ & $781$     & $1661$      & $2641$     & $3721$     \\
$2$     & $400$  & $1680$  & $1281$  & $7.29\,10^{-2}$ & $2961$    & $6321$      & $10081$    & $14241$    \\
$3$     & $900$  & $3720$  & $2821$  & $4.86\,10^{-2}$ & $6541$    & $13981$     & $22321$    & $31561$    \\
$4$     & $1600$ & $6560$  & $4961$  & $3.64\,10^{-2}$ & $11521$   & $24641$     & $39361$    & $55681$    \\
$5$     & $2500$ & $10200$ & $7701$  & $2.92\,10^{-2}$ & $17901$   & $38301$     & $61201$    & $--$       \\ 
$6$     & $3600$ & $14640$ & $11041$ & $2.43\,10^{-2}$ & $25681$   & $54961$     & $87841$    & $--$       \\ 
$7$     & $4900$ & $19880$ & $14981$ & $2.08\,10^{-2}$ & $34861$   & $74621$     & $119281$   & $--$       \\ 
$8$     & $6400$ & $25920$ & $19521$ & $1.82\,10^{-2}$ & $45441$   & $97281$     & $155521$   & $--$       \\ 
\end{tabular}
\end{center}
\end{small}
\end{table}
\begin{table}
\begin{small}
\caption{Geometric data and number of degrees of freedom of the sequence of randomized quadrilateral meshes.}
\label{table:mesh-data:randomized:quadrilateral}
\begin{center}
\begin{tabular}{
c|ccc|c|cccc}
$\ilev$ & $\nR$  & $\nF$   & $\nV$  & $\hmax$          & $\ndof{2}$ & $\ndof{3}$ & $\ndof{4}$ & $\ndof{5}$   \\
\hline
$0$     & $25$   & $60$    & $36$   & $3.311\,10^{-1}$ & $96$       & $216$      & $361$      & $531$      \\
$1$     & $100$  & $220$   & $121$  & $1.865\,10^{-1}$ & $341$      & $781$      & $1321$     & $1961$     \\
$2$     & $400$  & $840$   & $441$  & $9.412\,10^{-2}$ & $1281$     & $2961$     & $5041$     & $7521$     \\
$3$     & $900$  & $1860$  & $961$  & $6.130\,10^{-2}$ & $2821$     & $6541$     & $11161$    & $16681$    \\
$4$     & $1600$ & $3280$  & $1681$ & $4.693\,10^{-2}$ & $4961$     & $11521$    & $19681$    & $29441$    \\
$5$     & $2500$ & $5100$  & $2601$ & $3.808\,10^{-2}$ & $7701$     & $17901$    & $30601$    & $--$       \\ 
$6$     & $3600$ & $7320$  & $3721$ & $3.167\,10^{-2}$ & $11041$    & $25681$    & $43921$    & $--$       \\ 
$7$     & $4900$ & $9940$  & $5041$ & $2.751\,10^{-2}$ & $14981$    & $34861$    & $59641$    & $--$       \\ 
$8$     & $6400$ & $12960$ & $6561$ & $2.389\,10^{-2}$ & $19521$    & $45441$    & $77761$    & $--$       \\ 
\end{tabular}
\end{center}
\end{small}
\end{table}

\clearpage

\end{document}